\documentclass{article}
\usepackage{times}
\usepackage{empheq}
\usepackage{amsmath, amssymb}
\usepackage{amsthm}
\usepackage{bbm}
\usepackage[T1]{fontenc}
\usepackage[english]{babel}
\usepackage[utf8]{inputenc}
\usepackage{hyperref}
\usepackage{mathrsfs}
\usepackage{mathabx}
\usepackage{color}
\usepackage{anysize}
\usepackage{apacite}
\usepackage{natbib}
\marginsize{2.5cm}{2.5cm}{2cm}{2cm}
\usepackage[toc,page]{appendix}
\usepackage{url}
\usepackage{bm}
\usepackage{mathtools}
\usepackage[font={small,sc}]{caption}
\usepackage{ragged2e}
\usepackage[shortlabels]{enumitem}
\usepackage{algorithm}
\usepackage{algpseudocode}
\usepackage{algorithmicx}
\usepackage{comment}    
\usepackage{lineno}
\usepackage{tensor}
\usepackage{subcaption}

\usepackage{makecell}
\usepackage{multirow}
\usepackage{array}

\algnewcommand\algorithmicinput{\textbf{Input:}}
\algnewcommand\Input{\item[\algorithmicinput]}
\algnewcommand\algorithmicdefine{\textbf{Define}}
\algnewcommand\Define{\State\algorithmicdefine \ }
\algnewcommand\algorithmicverify{\textbf{Verify}}
\algnewcommand\Verify{\State\algorithmicverify \ }

\newtheorem{theorem}{Theorem}[section]
\newtheorem{definition}{Definition}[section]
\newtheorem{proposition}{Proposition}[section]

\newtheorem{lemma}{Lemma}[section]

\AtBeginDocument{}

\DeclareMathOperator\diam{diam}
\DeclareMathOperator\Span{span}

\newcommand{\Rd}{\mathbb{R}^{d}}
\newcommand{\Rm}{\mathbb{R}^{m}}

\newcommand{\RdxRd}{\Rd\times\Rd}
\newcommand{\RdxRm}{\Rd\times\Rm}
\newcommand{\RmxRm}{\Rm\times\Rm}

\newcommand{\Borel}[1]{\mathcal{B}(#1)}
\newcommand{\BoundedBorel}[1]{\mathcal{B}_{B}(#1)}

\newcommand{\BorelRd}{\mathcal{B}(\mathbb{R}^{d})}

\newcommand{\BoundedBorelRd}{\mathcal{B}_{B}(\mathbb{R}^{d})}

\newcommand{\Cov}{\mathbb{C}ov}
\newcommand{\Var}{\mathbb{V}ar}
\newcommand{\OmegaAP}{(\Omega , \mathcal{A} , \mathbb{P})}
\newcommand{\LtwoOmega}{L^{2}(\Omega , \mathcal{A} , \mathbb{P})}

\newcommand{\R}{\mathbb{R}}

\newcommand{\DRd}{\mathscr{D}(\Rd)}

\title{Karhunen-Loève expansion of random measures}

\author
{Ricardo Carrizo Vergara$^1$\\
     \normalsize{$^{1}$Université de Paris II Panthéon-Assas, Paris, France}\\
     \textit{racarriz@uc.cl}
}

\date{June 2025}

\begin{document}

\maketitle

\noindent {\bf Abstract} \quad  We present an orthogonal expansion for real, function-regulated, second-order random measures over $\Rd$ with measure covariance. Such a expansion, which can be seen as a Karhunen-Loève decomposition, consists in a series of deterministic real measures weighted by uncorrelated real random variables with the variances forming a convergent series. The convergence of the series is in a mean-square sense stochastically and against measurable bounded test functions (with compact support if the random measure is not finite) in the measure sense, which implies set-wise convergence. This is proven taking advantage of the extra requirement of having a covariance measure over $\RdxRd$ describing the covariance structure of the random measure, for which we also provide a series expansion. These results cover for instance the cases of Gaussian White Noise, Poisson and Cox point processes, and can be used to obtain expansions for trawl processes.
\bigskip


\section*{Introduction}
\label{Sec:Intro}

Karhunen-Loève (KL) expansions are an important tool for the analysis of stochastic processes, both in theory and practice. In a general non-rigorous manner, a KL expansion consists in a series representation for a random object $X$ taking values in a (real) vector space $E$, the representation being of the form
\begin{equation}
\label{Eq:GeneralKLX}
X = \sum_{n} X_{n}e_{n},
\end{equation}
where $(X_{n})_{n}$ is a collection of uncorrelated real random variables with $\sum_{n} \sigma_{X_{n}}^{2} < \infty $, $\sigma_{X_{n}}^{2} = \Var(X_{n})$, and $(e_{n})_{n}$ is a linearly independent collection of vectors in $E$. The case where $E$ is finite dimensional is commonly known as principal components analysis. When $E$ has infinite dimension, the most studied case is when $E$ is a separable Hilbert space, for which the vectors $(e_{n})_{n}$ form an orthonormal basis. The most basic scenario is when $X = (X(t))_{t \in [a,b]}$ is a mean-square continuous real stochastic process over a compact interval $[a,b] \subset \R$, in which $E = L^{2}([a,b])$ is used as basis Hilbert space which contains the continuous functions. The convergence of the series \eqref{Eq:GeneralKLX} has to be specified, both in the stochastic sense as a series of random objects, and in the sense of the space $E$ for which a topology must be made precise. The stochastic convergence of KL expansions is taken to be in mean-square. The case with $E$ Hilbert provides a direct topology on $E$ for the convergence (either the norm or weak topologies can be used). However, sometimes one can prove a stronger convergence than the one of the underlying Hilbert space. In the example of a mean-square continuous stochastic process over $[a,b]$, one uses $E = L^{2}([a,b])$, but Mercer's Theorem allows to conclude a stronger uniform-over-$[a,b]$-mean-square convergence. For general references on KL expansions and its applications, see \citep{loeve1978probability,wang2008karhunen,red2009elements}  .

This work focuses on KL expansions for random measures of a certain, very general kind. More precisely, Theorems \ref{Theo:KLMeasure} and \ref{Theo:KLExpansionMnonfinite} show that if $M$ is a (zero mean) \textit{measure-covariance second-order random measure} over $\Rd$ that is function-regulated (i.e., $M/f$ is finite for a given strictly positive function $f$, the case $f=1$ covering the case of finite measures), then $M$ admits the expansion
\begin{equation}
\label{Eq:MExpansionIntro}
M = \sum_{n \in \mathbb{N}}X_{n} \mu_{n},
\end{equation}
where $(X_{n})_{n \in \mathbb{N}}$ is a sequence of uncorrelated random variables with $\sum_{n}\sigma_{X_{n}}^{2} < \infty$, and, for the indexes $n$ such that $\sigma_{X_{n}}^{2} > 0 $, $(\mu_{n})_{n}$ form a collection of linearly independent real measures.\footnote{When $\sigma_{X_{n}}^{2} = 0$, the object $\mu_{n}$ may not be a measure but it does not really intervene in expansion \eqref{Eq:MExpansionIntro} since $X_{n}=0$ in such a case.} Theorem \ref{Theo:KLMeasure} covers the case of $M$ finite, which is where the real mathematical difficulty is present. Theorem \ref{Theo:KLExpansionMnonfinite} covers the function-regulated case, and it is essentially a corollary of Theorem \ref{Theo:KLMeasure}. The convergence of \eqref{Eq:MExpansionIntro} is in the sense
\begin{equation}
\label{Eq:ConvergenceIntro}
\mathbb{E}( |\langle M , \varphi \rangle - \sum_{j\leq n} X_{j} \langle \mu_{j} , \varphi \rangle |^{2} ) \xrightarrow[n \to \infty]{} 0, 
\end{equation} 
for every $\varphi$ measurable and bounded when $M$ is finite and every $\varphi$ measurable, bounded and compactly supported when $M$ is function-regulated. Here $\langle M , \varphi \rangle $ denotes the integral $\int_{\Rd}\varphi dM $.

The concept of a \textit{measure-covariance second-order random measure} deserves an explanation. First, we focus on a second-order random measure, which is a stochastic process indexed by bounded Borel sets, $M = (M(A))_{A \in \BoundedBorelRd}$ such that $M(A) \in \LtwoOmega$ for every $A$, with $\OmegaAP$ some probability space, and the application $A \mapsto M(A)$ is $\sigma$-additive. Note that this does not imply that $M$ is a \textit{measure-valued random variable}, that is, for a given $\omega \in \Omega$ the function $A \mapsto M(A)(\omega)$ is not necessarily a real measure over $\Rd$, nor almost surely in $\omega$. This measure-sample path definition is the one used by much of the current literature on random measures (see the introductory chapter in \citep{kallenberg2017random}), but it fails to cover very important cases such as Gaussian White Noise. Second-order random measures do contain Gaussian White Noise plus many other useful examples which we will mention further, but the literature on such random measures is more scarce; see  \citep{morando1969mesures,thornett1979class,rao2012random} as examples of general works using this concept. Now, the extra important adjective \textit{measure-covariance} comes from the very important assumption that there exists a measure $C_{M}$ over $\RdxRd$, called the \textit{covariance measure}, which satisfies
\begin{equation}
\label{Eq:CovMeasureMIntro}
\Cov( M(A) , M(B) ) = C_{M}(A \times B ), \quad \forall A,B \in \BoundedBorelRd.
\end{equation}
In general, a second-order random measure $M$ has its covariance structure determined by a \textit{bi-measure}, that is, the function $(A,B) \mapsto \Cov(M(A) , M(B))$ is a measure in one component when the other one is fixed. It is known \citep[Section 2.2, Example 2]{rao2012random} that a bi-measure is not generally identifiable with a measure over $\RdxRd$ as in \eqref{Eq:CovMeasureMIntro}. Therefore, the measure-covariance assumption is an extra regularity criteria which allows to obtain more conclusions. For instance, the existence of the total variation measure $|C_{M}|$ helps, as we shall see, to prove a semi-stochastic Fubini Theorem for random measures (Theorem \ref{Theo:StochasticFubini}, used mainly as an auxiliary result), and to prove the convergence type mode \eqref{Eq:ConvergenceIntro} thanks to the use of Lusin's Theorem. Assumption \eqref{Eq:CovMeasureMIntro} is still quite weak general and covers essentially every second-order random measure used in practice: we give some examples in Section \ref{Sec:ExamplesRandomMeasures}. We refer to \citep{borisov2006constructing,kruk2010malliavin} as other works where assumption \eqref{Eq:CovMeasureMIntro} is used. Expansion \eqref{Eq:MExpansionIntro} also implies an expansion for the covariance measure $C_{M}$, which is here specified in Proposition \ref{Prop:ExpansionCMabsVar}.

Orthogonal expansions for second-order random measures have been explored in particular cases. For Gaussian White Noise and other orthogonal random measures, orthogonal expansions can be obtained quite immediately (see Section \ref{Sec:WhiteNoise}). Note also that a second-order random measure $M$ can be interpreted as a generalized stochastic process (or random distribution, see \citep{gelfand1964generalized,ito1954stationary}) by focusing on the random variables $\langle M , \varphi \rangle$ for every $\varphi$ smooth with compact support. For such stochastic objects, orthogonal expansions such as \eqref{Eq:MExpansionIntro} are known. \citep{meidan1979reproducing} covers the case of a generalized stochastic process over a bounded subset of $\Rd$, the objects $(\mu_{n})_{n}$ being distributions in such case. \citep{carrizov2021generalized} explores the case of tempered random distributions over the whole space $\Rd$, in which second-order random measures regulated by polynomials are covered. The actual contribution of the present work is double: the demonstration that for a function-regulated measure-covariance second-order random measure $M$ the objects $(\mu_{n})_{n}$ in the KL expansion are \textit{measures} and not general distributions, all of them being regulated by the same function that regulates $M$; and the important convergence against measurable functions \eqref{Eq:ConvergenceIntro}, which is stronger that convergence against smooth functions and implies set-wise convergence. In what concerns KL expansions of more abstract stochastic objects with values in more general topological spaces than the classical Hilbert space case, we refer to \cite{bay2019karhunen} for separable Banach spaces, \cite{rajput1972gaussian} for separable Fréchet spaces, and \cite{peccati2010decompositions} for compact topological groups.

This work is organized as follows. In Section \ref{Sec:KLHilbert} we give the basis of KL expansions with respect to a Hilbert space. In Section \ref{Sec:RandomMeasures} we introduce random measures. Since the setting is not completely standard (use of $\delta$-rings and random measures in a particular sense), Sections \ref{Sec:DetMeasures} and \ref{Sec:DefRandomMeasures} introduce notations, basic notions and key properties of both deterministic and random measures over $\Rd$. Semi-stochastic Fubini Theorem \ref{Theo:StochasticFubini} is also here presented. Subsection \ref{Sec:ExamplesRandomMeasures} gives important examples of widely studied random measures for which a KL expansion as here presented can be obtained\footnote{We do not explicitly obtain their expansions here, we only mention them as examples covered by the results.} such as White Noise, orthogonal random measures, Poisson and Cox processes and the derivative of fractional Brownian motion with Hurst index $H \geq \frac{1}{2}$. In Section \ref{Sec:ExpansionRandomMeasures} we present the main Theorems \ref{Theo:KLMeasure} and \ref{Theo:KLExpansionMnonfinite}, Proposition \ref{Prop:ExpansionCMabsVar}, and their proofs. We end in Section \ref{Sec:Remarks} with some concluding remarks and comments about these results and ideas for future research. Namely, we discuss: the Hilbert space $E$ which contains the finite measures over $\Rd$ implicitly used as basis for the KL expansions; the uniqueness of these expansions; details in the Gaussian case; applicability for non-function-regulated random measures; and ideas for obtaining KL expansions for some non-mean-square-continuous stochastic processes over $\Rd$, such as trawl processes.

\textbf{Notations and conventions.} $\mathbf{1}_{A}$ denotes the indicator function of the set $A$. $\| \cdot \|_{\infty}$ denotes the supremum norm. $\DRd$ denotes the space of (real) smooth compactly supported test-functions over $\Rd$ typically used in Distribution Theory. The Lebesgue measure over $\Rd$ is denoted $\ell^{\otimes d}$. All random variables are supposed to be defined over a common probability space $\OmegaAP$. A stochastic process is understood as a family of random variables indexed by an arbitrary non-empty set. We do not make precise the laws of the random variables involved (the Gaussian case is a particular one which can be used as a reference example). Equality between random variables is always understood in an a.s. sense, and equality between stochastic processes is understood as one being a modification of the other.

\section{Karhunen-Loève Expansion}
\label{Sec:KLHilbert}

Let us give the details about KL expansions with respect to a Hilbert space. Let $E$ be a real separable Hilbert space, with inner-product $(\cdot , \cdot)_{E}$. Let $X : E \to \LtwoOmega$ be a linear and continuous real mapping satisfying that there exists an orthonormal basis $(e_{n})_{n \in \mathbb{N}} \subset E$ such that
\begin{equation}
\label{Eq:SeriesTrace}
\sum_{n \in \mathbb{N}} \mathbb{E}( |X(e_{n})|^{2} ) < \infty.
\end{equation}
If this holds, we say that $X$ has a \textit{traceable KL expansion with respect to $E$}. In such a case there exist an orthonormal basis of $E$, say $(f_{n})_{n \in \mathbb{N}}$, and a sequence of uncorrelated random variables $(X_{n})_{n \in \mathbb{N}}$, such that
\begin{equation}
\label{Eq:XeExpansion}
X(e) = \sum_{n \in \mathbb{N}}X_{n}\left( f_{n} , e \right)_{E}, \quad \forall e \in E,
\end{equation}
the convergence of the series being in a mean-square sense. Note that we have used the weak topology on $E$ for the convergence criterion. The vectors $(f_{n})_{n \in \mathbb{N}}$ are the eigenvectors of the \textit{covariance operator} induced by the covariance of $X$: if $K_{X} : E \times E \to \mathbb{R}$ is the covariance Kernel of $X$, that is
\begin{equation}
K_{X}(e,f) = \mathbb{E}(X(e)X(f)),
\end{equation}
then $K_{X}$ is bilinear, positive-semidefinite and continuous (since $X$ is continuous). By Riesz Representation, for every $e \in E$ there exists an element $Q_{X}(e) \in E $ such that
\begin{equation}
K_{X}(e,f) = \left( Q_{X}(e) , f \right)_{E}, \quad \forall f \in E.
\end{equation}
The so-induced operator $Q_{X} : E \to E $ is called the covariance operator of $X$. This operator is linear, continuous, positive-semidefinite, and by \eqref{Eq:SeriesTrace} it is also trace-class \citep[Theorem VI.18]{reed1980methods}. Hence, it has a spectral decomposition in an orthonormal basis of eigenvectors $(f_{n})_{n \in \mathbb{N}}$, with corresponding positive eigenvalues $(\sigma_{X_{n}}^{2})_{n \in \mathbb{N}}$ which form a convergent series \citep[Theorem VI.21]{reed1980methods}:
\begin{equation}
\sigma_{X_{n}}^{2}f_{n} = Q_{X}(f_{n}) \quad ; \quad \sum_{n \in \mathbb{N}} \sigma_{X_{n}}^{2} < \infty.
\end{equation}
The random variables $(X_{n})_{n \in \mathbb{N}}$ are given by $X_{n} := X(f_{n})$, for which we have $\Cov(X_{n},X_{m})= \sigma_{X_{n}}^{2} \delta_{n,m}$.

Let us study a particular example which we apply in this work. Let $(U(x))_{x \in \Rd}$ be a real mean-square continuous stochastic process over $\Rd$. Let $C_{U}(x,y) = \mathbb{E}(U(x)U(y))$ be its covariance function, which is continuous over $\RdxRd$. Let $\nu$ be a positive finite measure over $\Rd$ such that
\begin{equation}
\label{Eq:Cuxxdnufinite}
\int_{\Rd}C_{U}(x,x)d\nu(x) < \infty.
\end{equation}
From Cauchy-Schwarz inequality and the positive-semidefiniteness of $C_{U}$, \eqref{Eq:Cuxxdnufinite} implies
\begin{equation}
\label{Eq:CuMuMuIntegral}
\int_{\RdxRd} |C_{U}(x,y)||\varphi(x)||\phi(y)| d(\nu\otimes\nu)(x,y) < \infty, \quad \forall \varphi, \phi \in L^{2}(\Rd,\nu).
\end{equation}
It is known (see the details in Appendix \ref{App:FubiniProof}, use $\varphi d\nu$ as measure in  Lemma \ref{Lemma:IntZmuUnbounded}) that this condition allows to properly define the stochastic integrals
\begin{equation}
\label{Eq:DefUtilde}
\tilde{U}(\varphi) := \int_{\Rd} U(x) \varphi(x) d\nu(x), \quad \forall \varphi \in L^{2}(\Rd,\nu).
\end{equation}
Hence, one can re-define $U$ as a process indexed by functions in the separable Hilbert space $L^{2}(\Rd,\nu)$. The so-defined application $\tilde{U} : L^{2}(\Rd,\nu) \to \LtwoOmega$ is continuous. The covariance operator $Q_{\tilde{U}}$ is given by
\begin{equation}
\label{Eq:DefCovOperaturUL2nu}
Q_{\tilde{U}}(\varphi) = \int_{\Rd}C_{U}(\cdot , y)\varphi(y)d\nu(y),
\end{equation}
which by \eqref{Eq:Cuxxdnufinite} is trace-class \citep{brislawn1991traceable}. $\tilde{U}$ has then a traceable KL expansion with respect to $L^{2}(\Rd,\nu)$:
\begin{equation}
\label{Eq:ExpansionUL2nu}
\int_{\Rd}U(x)\varphi(x)d\nu(x) = \sum_{n \in \mathbb{N}} X_{n}( f_{n} , \varphi  )_{L^{2}(\Rd,\nu)}, \quad \forall \varphi \in L^{2}(\Rd,\nu),
\end{equation} 
with $(f_{n})_{n \in \mathbb{N}}$ the orthonormal basis of eigenfunctions of $Q_{\tilde{U}}$ and
\begin{equation}
\label{Eq:DefXnUfnnu}
X_{n} = \int_{\Rd}U(x)f_{n}(x)d\nu(x), \quad \forall n \in \mathbb{N}.
\end{equation}

\section{Random measures}
\label{Sec:RandomMeasures}

The proofs of the claims exposed in this Section are present in Appendix \ref{App:ProofsSectionRM}.

\subsection{Reminders on measures over $\Rd$ and their anti-derivatives}
\label{Sec:DetMeasures}

We denote $\BorelRd$ the Borel $\sigma$-algebra of $\Rd$ and $\BoundedBorelRd$ the $\delta$-ring of bounded Borel subsets of $\Rd$. By a  \textit{measure over} $\Rd$, we mean a \textit{real} application $\mu : \BoundedBorelRd \to \R$ which is $\sigma$-additive over $\BoundedBorelRd$. This implies that $\mu$ is locally-finite, but $\mu$ may not be defined over unbounded sets (some authors use the term \textit{pre-measure} for this object \citep{kupka1978caratheodory}). $\mu$ is called positive if it takes only non-negative values. The total-variation measure of $\mu$, noted $|\mu|$, is the smallest positive measure such that $|\mu(A)| \leq |\mu|(A) $ for all $A \in \BoundedBorelRd$ \citep[Chapter 6]{rudin1987real}. If $|\mu|$ can be extended finitely and $\sigma$-additively to $\BorelRd$ (hence $|\mu|(\Rd) < \infty $), then $\mu$ is said to be finite, and it can be extended uniquely and $\sigma$-additively to $\BorelRd$. The space of measures (resp. finite measures) over $\Rd$ is denoted $\mathscr{M}(\Rd)$ (resp. $\mathscr{M}_{F}(\Rd)$). The space of (real Borel) measurable functions over $\Rd$ is denoted $\mathcal{M}(\Rd)$.
$\mathcal{M}_{B}(\Rd)$ and $\mathcal{M}_{B,c}(\Rd)$ denote the subspaces of $\mathcal{M}(\Rd)$ consistent of bounded and bounded compactly supported functions respectively. A function $f \in \mathcal{M}(\Rd)$ is said to be integrable with respect to $\mu \in \mathscr{M}(\Rd)$ if $|f|$ is Lebesgue integrable with respect to $|\mu|$. In such case, we note $\langle \mu , f \rangle := \int_{\Rd}f d\mu = \int_{\Rd}f(x)d\mu(x).$ We remark that the total variation measure $|\mu|$ can be expressed as
\begin{equation}
\label{Eq:TotalVariationAlternative}
|\mu|(A) = \sup_{\varphi \in \mathcal{M}(\Rd) , |\varphi | = \mathbf{1}_{A}} | \langle \mu , \varphi \rangle |, \quad \forall A \in \BoundedBorelRd.
\end{equation}

We recall the useful Lusin's Theorem, considered in a simplified version over $\Rd$ \citep[Theorem 7.10]{folland1999real}: 
\begin{theorem}[\textbf{Lusin}]
\label{Theo:Lusin}
Let $\mu \in \mathscr{M}_{F}(\Rd)$ and $\psi \in \mathcal{M}(\Rd)$. Then, for every $\epsilon > 0 $ there exists a closed set $E \subset \Rd$ such that $\psi$ is continuous over $E$ (with the subspace topology) and $|\mu|(E^{c}) < \epsilon $.
\end{theorem}

One special property of measures over $\Rd$ is that they are derivatives in distributional sense of regular functions. Moreover, if the measure is finite those primitives grow in a controlled manner. Consider thus the following (double) anti-derivative operator $\mathcal{O} : \mathscr{M}_{F}(\Rd) \to C(\Rd)$:
\begin{equation}
\label{Eq:DefO}
\mathcal{O}( \mu )(\vec{x}) := \int^{\vec{x}}_{0} \mu( ( - \infty , \vec{u} ] ) d\vec{u}, \quad \forall \vec{x} \in \Rd, 
\end{equation}
where we have used the abbreviated notations
\begin{equation}
\label{Eq:IntegralVecXnotation}
\int^{\vec{x}}_{0} ( \cdot ) \ d\vec{u} := \int_{0}^{x_{1}} \int_{0}^{x_{2}} ... \int_{0}^{x_{d}} ( \cdot ) \ du_{d}... du_{2} du_{1} \quad ; \quad ( - \infty , \vec{x} ] := (-\infty , x_{1} ] \times (-\infty , x_{2} ] \times ... \times (-\infty , x_{d} ], 
\end{equation}
for every $\vec{x} = (x_{1} , ... , x_{d} ) \in \Rd$. Note that the function $\vec{u} \mapsto \mu( ( - \infty , \vec{u} ] )$ is bounded (since $\mu$ is finite) and càdlàg in each component when the others are fixed, therefore the iterated integrals in \eqref{Eq:DefO} are simple Riemann integrals and thus $\mathcal{O}(\mu)$ is a continuous function. The following bound holds for $\mathcal{O}(\mu)$:

\begin{equation}
\label{Eq:OmuBounded}
|\mathcal{O}(\mu)(\vec{x})| \leq |x_{1}|...|x_{d}| |\mu|(\Rd), \quad \forall \vec{x} = (x_{1} , ... , x_{d}) \in \Rd.
\end{equation}
$\mathcal{O}$ is an anti-derivative operator in the sense that $\frac{\partial^{2d} \mathcal{O}(\mu)  }{
\partial x_{1}^{2} ... \partial x_{d}^{2}} = \mu$ in distributional sense over $\Rd$, that is
\begin{equation}
\label{Eq:DerOmuMuIntegral}
\int_{\Rd}  \mathcal{O}(\mu)(x)   \frac{\partial^{2d} \varphi  }{
\partial x_{1}^{2} ... \partial x_{d}^{2}}(x)dx = \int_{\Rd} \varphi(x) d\mu(x), \quad \forall \varphi \in \DRd.
\end{equation}

\subsection{Measure-covariance random measures and properties}
\label{Sec:DefRandomMeasures}

\begin{definition}
\label{Def:RM}
A centred measure-covariance second-order random measure (from now on \textit{m-cov random measure}) over $\Rd$ is a zero-mean stochastic process indexed by the bounded Borel sets $M = (M(A))_{A \in \BoundedBorelRd}$ such that there exists $C_{M} \in \mathscr{M}(\RdxRd)$ such that
\begin{equation}
\label{Eq:DefCM(AxB)}
\mathbb{E}(M(A)M(B) ) = C_{M}(A\times B), \quad \forall A,B \in \BoundedBorelRd.
\end{equation}
\end{definition}

The first implication of Definition \ref{Def:RM} is the following.
\begin{proposition}
\label{Prop:MsigmaAdditive}
$M$ is a $\sigma$-additive function from $\BoundedBorelRd$ to $\LtwoOmega$.
\end{proposition}

In other words, $M$ is an $\LtwoOmega$-valued (locally finite) measure over $\Rd$. The extra adjective \textit{m-cov} is added because of the identification of the covariance of $M$ to the \textit{covariance measure} $C_{M}$, which, as mentioned in the introduction, does not apply for a general second-order random measure \citep[Chapter 2, Example 2]{rao2012random}. It is clear that covariance measures are symmetric in the sense $C_{M}(A \times B ) = C_{M}( B \times A )$. It is possible to verify that $|C_{M}|$ is also a symmetric measure.  Covariance measures are positive-semidefinite in the sense
\begin{equation}
\label{Eq:CMdefpos}
\langle C_{M} , \varphi \otimes \varphi \rangle \geq 0 , \quad \forall \varphi \in \mathcal{M}_{B,c}(\Rd).
\end{equation}
Conversely, every symmetric measure over $\RdxRd$ satisfying \eqref{Eq:CMdefpos} is the covariance measure of an m-cov random measure.\footnote{Construct a Gaussian m-cov random measure using Kolmogorov Extension Theorem.} If $\varphi \in \mathcal{M}(\Rd)$ is such that
\begin{equation}
\label{Eq:VarphiIntegrableM}
\langle |C_{M}| , |\varphi|\otimes |\varphi| \rangle < \infty,
\end{equation}
then the stochastic integral
\begin{equation}
\label{Eq:IntegralMVarphi}
\langle M , \varphi \rangle := \int_{\Rd}\varphi(x) dM(x)
\end{equation}
can be uniquely defined as a random variable in $\LtwoOmega$. This is just an example of the \textit{Dunford-Schwartz} integral of $\varphi$ with respect to the $\LtwoOmega$-valued measure $M$; see \citep[Chapter 2]{rao2012random} for an effective introduction, \citep[Section IV.10]{dunford1958linear} for the details, and \citep[Proposition 3.3.1]{carrizov2018development} for the sufficiency of condition \eqref{Eq:VarphiIntegrableM}. If $\varphi$ and $\phi$ satisfy \eqref{Eq:VarphiIntegrableM}, then
\begin{equation}
\label{Eq:CovarianceIntegrals}
\mathbb{E}\left(  \langle M , \varphi \rangle \langle M , \phi \rangle \right) = \langle C_{M} , \varphi \otimes \phi \rangle.
\end{equation}

The next theorem, which will play an auxiliary role, is called here \textit{semi-stochastic Fubini Theorem}, since it provides sufficient conditions under which we can switch integral signs when one of the integrating measures is random and the other is not. Other stochastic Fubini theorems can be found in the literature but usually with diverse sample path, predictability or martingale-type conditions (\citep[Theorem 7.4.10]{rao2012random}, \citep[Section 4.5]{daprato2014stochastic}, \citep{veraar2012stochastic}), which are not the focus here. The version here presented aims to provide conditions on $C_{M}$ so stochastic integrals can be defined with Riemann sums, without extra requirements on the sample paths of $M$.

\begin{theorem}[\textbf{Semi-stochastic Fubini}]
\label{Theo:StochasticFubini}
Let $M$ be an m-cov random measure over $\Rd$ with covariance measure $C_{M}$ and let $\mu \in \mathscr{M}(\Rm)$. Let $\psi \in \mathcal{M}(\RdxRm)$ such that
\begin{enumerate}[\itshape(i)]
\item  $  \int_{\RdxRd \times \Rm\times\Rm }   |\psi(x,u)| |\psi(y,v) | d |C_{M}|\otimes |\mu|\otimes |\mu|  (x,y,u,v)  < \infty.  $
\label{It:Fubini1}
\item  The function $ (u,v) \mapsto  \langle \ |C_{M}| \ , \ |\psi|( \cdot  , u  ) \otimes |\psi|( \cdot , v  ) \ \rangle $ is locally bounded and there exists $E \in \Borel{\Rm}$ with $|\mu|(E^{c})=0$ such that the function $ (u,v) \mapsto  \langle \ C_{M} \ , \ \psi( \cdot  , u  ) \otimes \psi( \cdot , v  ) \ \rangle $ is continuous over $E\times E$. \label{It:Fubini2}
\end{enumerate} 
Then, 
\begin{equation}
\label{Eq:FubiniTheorem}
\int_{\Rd}\int_{\Rm} \psi(x,u) d\mu(u) dM(x) = \int_{\Rm}\int_{\Rd} \psi(x,u) dM(x) d\mu(u) . 
\end{equation}
\end{theorem}

Let us now focus on the case of finite random measures.

\begin{definition}
\label{Def:FiniteRM}
An m-cov random measure $M$ over $\Rd$ is said to be finite if its covariance measure $C_{M}$ is finite.
\end{definition}

When $M$ is finite, its definition can be extended uniquely, finitely and $\sigma$-additively to the whole Borel $\sigma$-algebra $\BorelRd$, the random variable $M(\Rd)$ having finite variance. There is also an extra regularity property which holds for finite random measures.

\begin{proposition}
\label{Prop:CMprimitiveContinuous}
Let $M$ be an m-cov finite random measure over $\Rd$. Then, the function over $\Rd\times\Rd$
\begin{equation}
\label{Eq:CovPrimitiveCM}
(\vec{u} , \vec{v} ) \mapsto C_{M}\left( (-\infty , \vec{u} ] \times (-\infty , \vec{v}] \right)
\end{equation}
is continuous over a set of the form $E\times E$, with $E \in \BorelRd$ such that $\ell^{\otimes d}(E^{c}) = 0$.
\end{proposition}

Let us now define the application of the anti-derivative operator $\mathcal{O}$ to a finite m-cov random measure $M$. The application $\vec{u} \mapsto M( (-\infty , \vec{u} ]) $ defines a stochastic process over $\Rd$ whose covariance function is \eqref{Eq:CovPrimitiveCM}, being thus mean-square continuous outside a set of null Lebesgue measure and with bounded covariance. Thus, the stochastic integral (see Lemma \ref{Lemma:IntZmuRiemann})
\begin{equation}
\label{Eq:DefOM}
\mathcal{O}( M )(\vec{x}) := \int_{0}^{\vec{x}} M( (-\infty , \vec{u}] ) d\vec{u}
\end{equation}
is well-defined through Riemann-alike approximations. This process has covariance function
\begin{equation}
\label{Eq:COM}
C_{\mathcal{O}( M )}(\vec{x} , \vec{y} ) = \int_{0}^{\vec{x}}\int_{0}^{\vec{y}} C_{M}( (-\infty , \vec{u} ]\times (-\infty , \vec{v}] ) d\vec{v} d\vec{u},
\end{equation}
which is a continuous function over $\RdxRd$ (it is actually the function $\mathcal{O}\otimes \mathcal{O}( C_{M} )$), therefore $\mathcal{O}(M)$ is mean-square continuous. In addition one has the bound
\begin{equation}
\label{Eq:COMBound}
|C_{\mathcal{O}( M )}(\vec{x} , \vec{y} )| \leq |x_{1}|...|x_{d}| |y_{1}| ... |y_{d}| |C_{M}|(\RdxRd), \quad \forall \vec{x} , \vec{y} \in \Rd.
\end{equation}
Finally, an application of semi-stochastic Fubini Theorem \ref{Theo:StochasticFubini} allows to conclude $\frac{\partial^{2d} \mathcal{O}(M)  }{
\partial x_{1}^{2} ... \partial x_{d}^{2}} = M$ in distributional sense over $\Rd$, that is, we have the equality between the stochastic integrals
\begin{equation}
\label{Eq:OMDerMIntegrals}
\int_{\Rd} \mathcal{O}(M)(x)  \frac{\partial^{2d} \varphi  }{
\partial x_{1}^{2} ... \partial x_{d}^{2}}(x) dx = \int_{\Rd} \varphi(x) dM(x), \quad \forall \varphi \in \DRd.
\end{equation}

\subsection{Examples of m-cov random measures}
\label{Sec:ExamplesRandomMeasures}

We provide some examples of commonly used random measures for which a KL expansion as it is presented here can be obtained. 

\subsubsection{White Noise and other orthogonal random measures}
\label{Sec:WhiteNoise}

A (non-necessarily Gaussian) White Noise over $\Rd$ is a centred m-cov random measure $W=(W(A))_{A \in \BoundedBorelRd}$ with covariance given by
\begin{equation}
\label{Eq:CovW}
\Cov( W(A) , W(B) ) = \ell^{\otimes d}( A \cap B ).
\end{equation}
The covariance measure of $W$ satisfies $ \langle C_{W} , \psi \rangle = \int_{\Rd} \psi(x,x) dx$ for every $\psi \in \mathcal{M}_{B,c}(\RdxRd)$. $C_{W}$  is a measure concentrated on the hyperplane $\{ x = y \} = \lbrace (x,y) \in \RdxRd \ \mid \ x = y \rbrace$, sometimes denoted $\delta(x-y)$. White Noise is a particular case of an orthogonal random measure. An orthogonal random measure is a centred m-cov random measure $M = (M(A))_{A \in \BoundedBorelRd}$ such that there exists  $\nu \in \mathscr{M}(\Rd)$ positive such that
\begin{equation}
\label{Eq:CovMOrthogonal}
\Cov(M(A) , M(B) ) = \nu(A\cap B).
\end{equation}
$C_{M}$ is also concentrated on $\{ x = y \}$ but with another weighting measure, having $\langle C_{M} , \psi \rangle = \int_{\Rd}\psi(x,x)d\nu(x)$. We denote in such case $C_{M}=\nu \delta(x-y)$. Orthogonal random measures have the property of assigning null covariance when evaluated over disjoint sets, or when integrated against functions with disjoint support. These kinds of random measures appear in the spectral analysis of stationary random random fields \citep{yaglom1987correlation}. In the stronger case where $M$ takes independent values at disjoint sets, $M$ is sometimes called a \textit{completely random measure} \citep{kingman1967completely,collet2021completely}, or an \textit{independently scattered random measure} \citep{passeggeri2020extension}. Lévy processes \citep{ken1999levy} can be seen as primitives of completely random measures (the \textit{Lévy basis}), and therefore their derivatives in distributional sense are orthogonal random measures in the sense \eqref{Eq:CovMOrthogonal} if the increments of the Lévy process have finite variance. KL expansions for Lévy processes have been worked out for example in \citep{hackmann2018karhunen}.

Orthogonal expansions for an orthogonal random measure $M$ can be obtained with relative ease. Let $(f_{n})_{n}$ to be an orthonormal basis of the space $L^{2}(\Rd,\nu)$. Then,
\begin{equation}
\label{Eq:ExpansionMOrthogonal}
\langle M , \varphi \rangle = \sum_{n \in \mathbb{N}} \epsilon_{n} ( f_{n} , \varphi  )_{L^{2}(\Rd , \nu)} , \quad \forall \varphi \in L^{2}(\Rd,\nu), 
\end{equation}
with $\epsilon_{n} = \langle M , f_{n} \rangle$. Note that in this case $\Var(\epsilon_{n}) = 1$, therefore the expansion is not traceable contrarily to the case of Section \ref{Sec:KLHilbert}. This can be arranged, for example, by multiplying each $\epsilon_{n}$ by a coefficient $\sigma_{n}>0$, with $\sum_{n} \sigma_{n}^{2} < \infty$, and then take $\mu_{n} = f_{n}/\sigma_{n}$ as functions in the expansion. Note that in such a case, $(\mu_{n})_{n}$ is not an orthonormal system of $L^{2}(\Rd,\nu)$ but of another more abstract Hilbert space, with respect to which $M$ has a traceable KL expansion (see further in Section \ref{Sec:HilbertApproach}). In order to identify expansion \eqref{Eq:ExpansionMOrthogonal} as a KL expansion such as the here developed, $f_{n}$ and $\mu_{n}$ must be interpreted as measures, not as functions.

Some orthogonal random measures provide the crucial example of $\LtwoOmega$-valued random measures that cannot be seen as random measures in the sense of random variables taking values in a space of measures or almost surely so. Over $\Rd$, independently scattered  measure-valued random variables must necessarily be a point process \citep{kingman1967completely}. In consequence, if $M$ is a Gaussian orthogonal random measure such that the weighting measure $\nu$ is not purely a  discrete measure, the sample paths of Gaussian orthogonal random measures have  almost surely unbounded variation \citep{horowitz1986gaussian}. This includes the case of Gaussian White Noise, as it is widely known \citep[Exercice 2.17]{oksendal2003stochastic}.

\subsubsection{Poisson and Cox point processes}
\label{Sec:PoissonCox}

A point process \citep{daley2006introduction} is a stochastic process indexed by the bounded Borel sets $(M(A))_{A \in \BoundedBorelRd}$ which can be represented as
\begin{equation}
\label{Def:PointProcess}
M(A) = \sum_{j \in \mathbb{N}} \delta_{X_{j}}(A),
\end{equation}
where $(X_{j})_{j \in \mathbb{N}}$ is a family of $\Rd$-valued random variables such that $M(A) < \infty$ almost surely. $M$ is called an inhomogeneous Poisson process if for every disjoint collection of bounded Borel sets $(A_{k})_{k}$, the random variables $(M(A_{k}))_{k}$ are independent Poisson random variables with $\mathbb{E}\left( M(A_{k}) \right) = \nu(A_{k})$ for some positive measure $\nu \in \mathscr{M}(\Rd)$ (the \textit{intensity} measure). From the independence at disjoint sets condition, the covariance structure of an inhomogeneous Poisson process is given by \eqref{Eq:CovMOrthogonal}, and thus $M - \nu$ (that is, centering $M$) is an orthogonal random measure. In consequence, orthogonal expansions of the form \eqref{Eq:ExpansionMOrthogonal} also hold for it, with $(f_{n})_{n}$ interpreted as measures.

Now, let $\Lambda : \BoundedBorelRd \mapsto \LtwoOmega$ be a \textit{positive} second-order random measure over $\Rd$, that is, $\Lambda(A) \geq 0 $ for every $A \in \BoundedBorelRd$. It is known \citep[Proposition 2.4]{rajput1989spectral} that in such case, the covariance of $\Lambda$ is always identified with a covariance measure $C_{\Lambda} \in \mathscr{M}(\RdxRd)$. In addition, the $\sigma$-additivity implies that the mean $\nu(A) := \mathbb{E}( \Lambda(A) )$ defines a measure $\nu \in \mathscr{M}(\Rd)$. Now, define $M$ such that, conditioned on $\Lambda$, $M$ is a Poisson point process with intensity $\Lambda$. Then, $M$ is another form of point process, commonly used in applications, called the Cox process \citep{cox1955some}. In such case $M - \nu$ is also an m-cov random measure, with covariance measure
\begin{equation}
\label{Eq:CovMCox}
C_{M} = \nu \delta(x-y)  + C_{\Lambda}.  
\end{equation}
Thus, $M$ has a richer covariance structure than a Poisson process, with an orthogonal random measure part $\nu \delta(x-y)$ plus an extra positive covariance $C_{\Lambda}$. The most popular Cox process among applications is the log-Gaussian Cox  process \citep{moller1998log}, where the random intensity is given by 
\begin{equation}
\label{Eq:LambdaLogGaussian}
\Lambda(A) = \int_{A} e^{Z(x)}dx,
\end{equation}
where $Z$ is some mean-square continuous Gaussian process. Note that random measures constructed from the integrals of an enough regular stochastic process with respect to a deterministic measure such as in \eqref{Eq:LambdaLogGaussian} also provide an example of m-cov random measures, see Appendix \ref{App:FubiniProof}.

\subsubsection{Derivatives of fractional Brownian motion}
\label{Sec:FracBrown}

Let $(B_{H}(t))_{t \geq 0}$ be a zero-mean $\R$-valued Gaussian process with covariance function
\begin{equation}
C_{B_{H}}(t,s) = \Cov(B_{H}(t) , B_{H}(s) ) = \frac{t^{2H} + s^{2H} - |t-s|^{2H}}{2},
\end{equation}
where $H \in (0,1)$. Then $B_{H}$ is called a \textit{fractional Brownian motion} and $H$ is called the Hurst index. If $H = \frac{1}{2}$, $B_{H}$ is a standard Brownian motion. Consider the case $H > \frac{1}{2}$. Let $\frac{d}{dt}B_{H}$ be the distributional derivative of $B_{H}$, whose covariance is given by \citep{borisov2006constructing}
\begin{equation}
\frac{\partial^{2}}{\partial t \partial s} C_{B_{H}} = H(2H-1)|t-s|^{2H - 2}, 
\end{equation}
which is not a continuous function but it is integrable over $[0,T]\times[0,T]$ for every $T > 0$. It follows that the covariance of $\frac{d}{dt}B_{H}$ can be identified with the measure
\begin{equation}
C_{\frac{d}{dt}B_{H}}(E) = H(2H-1)\int_{E}\frac{d(x,y)}{|x-y|^{2-2H}}, \quad \forall E \in \BoundedBorel{[0,\infty)\times[0,\infty)}. 
\end{equation}
$\frac{d}{dt}B_{H}$ is thus another example of an m-cov random measure. Note that some authors call $\frac{d}{dt}B_{H}$ a long-range dependence process \citep{gay1990class,anh1999possible}. It is known that the case $H > \frac{1}{2}$ is regular enough to develop an stochastic calculus around $B_{H}$ without requiring specialized techniques, contrarily to the Brownian motion case \citep{zahle1998integration}.

\section{Expansion of random measures}
\label{Sec:ExpansionRandomMeasures}

Now that every required definition and basic result is established, we present the KL expansion for finite random measures, which is the main result of this work.

\begin{theorem}[\textbf{Karhunen-Loève expansion of finite random measures}]
\label{Theo:KLMeasure} 
Let $M$ be an m-cov finite random measure over $\Rd$. Then, there exists a sequence of pairwise uncorrelated random variables with summable variances $(X_{n})_{n \in \mathbb{N}}$, and a linearly independent sequence of finite measures over $\Rd$, $(\mu_{n})_{n \in \mathbb{N}}$ such that
\begin{equation}
\label{Eq:MvarphiExpansion}
\langle M , \varphi \rangle = \sum_{n \in \mathbb{N}}X_{n}\langle \mu_{n} , \varphi \rangle, \quad \forall \varphi \in \mathcal{M}_{B}(\Rd),
\end{equation}
with the series being considered in a mean-square sense. 
\end{theorem}

The arguments behind the proof of Theorem \ref{Theo:KLMeasure} are actually simple. We first apply the anti-derivative operator $\mathcal{O}$ to $M$ in order to obtain an enough regular process for which a KL expansion with respect to some Hilbert space exists. Then, we derive it to retrieve $M$. This logic has been applied for the case of general tempered random distributions \citep{carrizov2021generalized}. The particularity here is the measure structure of the objects $(\mu_{n})_{n}$ and the convergence mode \eqref{Eq:MvarphiExpansion}, which requires extra attention. The proof will be split into a few Lemmas. New notations will be introduced and kept along the Lemmas. The reader may recognize very similar arguments to the proof of the classical KL expansion for mean-square continuous stochastic process over compact intervals \citep[Section 37.5]{loeve1978probability}.

\begin{lemma}
\label{Lemma:Nu}
There exists $\nu \in \mathscr{M}_{F}(\Rd)$ such that for every m-cov finite random measure $M$ over $\Rd$ the process $\mathcal{O}(M)$ has a KL expansion with respect to $L^{2}\left( \Rd , \nu \right).$
\end{lemma}

\textbf{Proof of Lemma \ref{Lemma:Nu}:} 
Consider the polynomial function $p: \Rd \to \R^{+}$ given by $p(\vec{x}) = \prod_{j=1}^{d}(1+|x_{j}|^{2})^{2}$. Consider the finite measure over $\Rd$
\begin{equation}
\label{Eq:Defnu}
d\nu(x) := \frac{dx}{p(x)}.
\end{equation}
Let $C_{\mathcal{O}(M)}$ be the covariance function of $\mathcal{O}(M)$. From bound \eqref{Eq:COMBound} we conclude  
\begin{equation}
\label{Eq:COMtraceIntegrableNu}
\int_{\RdxRd} C_{\mathcal{O}(M)}(x,x)d\nu(x) \leq |C_{M}|(\RdxRd)\left( \int_{\R}  \frac{t^{2}}{(1+t^{2})^{2}} dt\right)^{d} < \infty.
\end{equation}
$\mathcal{O}(M)$ has thus a traceable KL expansion with respect to $L^{2}(\Rd,\nu)$ (Section \ref{Sec:KLHilbert}), having thus
\begin{equation}
\label{Eq:OMexpansionL2nu}
\int_{\Rd}\mathcal{O}(M)(x) \varphi(x) d\nu(x) = \sum_{j \in \mathbb{N}} X_{j} ( f_{j} , \varphi )_{L^{2}(\Rd,\nu)} = \sum_{j \in \mathbb{N}} X_{j} \langle \nu , f_{j} \varphi \rangle, \quad \forall \varphi \in L^{2}(\Rd,\nu),
\end{equation}
being $(f_{j})_{j \in \mathbb{N}}$ the orthonormal basis of $L^{2}(\Rd,\nu)$ given by the eigenfunctions of the covariance operator of $\mathcal{O}(M)$, and $(X_{j})_{j\in \mathbb{N}}$ the associated uncorrelated random variables with variances $(\sigma_{X_{j}}^{2})_{j\in \mathbb{N}}$, satisfying
\begin{equation}
\label{Eq:ObjectsExpansionOML2nu}
\sigma_{X_{j}}^{2}f_{j} = \int_{\Rd}C_{\mathcal{O}(M)}( \cdot , y ) f_{j}(y) d\nu(y) \quad ; \quad X_{j} = \int_{\Rd}\mathcal{O}(M)(x)f_{j}(x)d\nu(x) \quad ; \quad \sum_{j \in \mathbb{N}} \sigma_{X_{j}}^{2} < \infty.\quad \blacksquare
\end{equation}

The measure $\nu$ in Lemma \ref{Lemma:Nu} is far from being unique: one can take any measure with fast-enough decreasing density so second-order primitives (in each component) of finite measures are integrable with respect to it (a Gaussian density works, for instance).  

Let us now fix $M$ as a given finite m-cov random measure.

\begin{lemma}
\label{Lemma:DerfjMeasures}
For every $j$ such that $\sigma_{X_{j}}^{2} > 0 $, the distribution $ \frac{\partial^{2d} f_{j}}{\partial x_{1}^{2} \ldots \partial x_{d}^{2}} $ is in $\mathscr{M}_{F}(\Rd)$. 
\end{lemma}

\textbf{Proof of Lemma \ref{Lemma:DerfjMeasures}:}
An arbitrary $f \in L^{2}(\Rd,\nu)$ determines a distribution over $\Rd$ through the application $\varphi \mapsto \int_{\Rd}f(x)\varphi(x)dx =  \langle fp\nu , \varphi \rangle$. Therefore, the derivatives $ \frac{\partial^{2d} f_{j}}{\partial x_{1}^{2} \ldots \partial x_{d}^{2}} $ are well defined as distributions over $\Rd$. The eigenvalue-eigenfunction relation implies for $\sigma_{X_{j}}^{2} > 0 $

\begin{equation}
\label{Eq:fjDevelop}
\begin{aligned}
f_{j}(\vec{x}) &= \frac{1}{\sigma_{X_{j}}^{2}} \int_{\Rd} C_{\mathcal{O}(M)}(\vec{x} , \vec{y} ) f_{j}(\vec{y})d\nu(\vec{y}) \\
&=  \frac{1}{\sigma_{X_{j}}^{2}} \int_{\Rd} \int_{0}^{\vec{x}} \int_{0}^{\vec{y}} C_{M}( (-\infty , \vec{u} ]\times(-\infty , \vec{v} ]  )d\vec{v} d\vec{u} f_{j}(\vec{y})d\nu(\vec{y}) \\
&= \frac{1}{\sigma_{X_{j}}^{2}} \int_{0}^{\vec{x}} \int_{\Rd} \int_{0}^{\vec{y}} C_{M}( (-\infty , \vec{u} ]\times(-\infty , \vec{v} ]  )d\vec{v}  f_{j}(\vec{y})d\nu(\vec{y}) d\vec{u} \\
&= \frac{1}{\sigma_{X_{j}}^{2}} \int_{0}^{\vec{x}} \int_{\Rd}  \mathcal{O}\left( C_{M}\left(  \left(-\infty , \vec{u} \right] \times \ \cdot \  \right) \right) (\vec{y}) f_{j}(\vec{y})d\nu(\vec{y}) d\vec{u},
\end{aligned}
\end{equation}
where we used (deterministic) Fubini Theorem\footnote{The classical Fubini Theorem for positive measures can be extended easily to the case of real measures over Euclidean spaces provided that the corresponding integrals using the total-variation of the measures involved are finite.} for changing the order of integration, and  $C_{M}\left( A \times  \ \cdot \    \right)$ stands for the measure $B \mapsto C_{M}\left(  A \times  B   \right)$ for any $A \in \BorelRd$. Inspired by this, we define
\begin{equation}
\label{Eq:Defmuj}
\mu_{j}(A) := \frac{1}{\sigma_{X_{j}}^{2}}\int_{\Rd} \mathcal{O}\left( C_{M}\left( A\times \ \cdot \  \right) \right)(\vec{y})  f_{j}(\vec{y})d\nu(\vec{y}), \quad \forall A \in \BorelRd.
\end{equation}
Given the property \eqref{Eq:OmuBounded} of the operator $\mathcal{O}$, the function $\mathcal{O}\big( C_{M}(A\times \cdot ) \big)$ is in $L^{2}(\Rd, \nu)$ and the integral \eqref{Eq:Defmuj} is thus well-defined. Since $\mathcal{O}$ is linear, the application $A \mapsto \mu_{j}(A)$ is additive. From bound \eqref{Eq:OmuBounded} we have
\begin{equation}
\label{Eq:BoundMuj}
|\mu_{j}(A)| \leq \frac{1}{\sigma_{X_{j}}^{2}} \underbrace{\int_{\Rd}|y_{1}|...|y_{d}| |f_{j}(\vec{y})|d\nu(\vec{y})}_{< \infty \hbox{ since } \vec{y} \mapsto |y_{1}|...|y_{d}|\  \in \ L^{2}(\Rd,\nu)} |C_{M}|(A \times \Rd).
\end{equation}
Since $|C_{M}|$ is a finite measure, if we take any sequence of Borel sets $(A_{n})_{n \in \mathbb{N}}$ such that $A_{n}\searrow \emptyset $ we have $|C_{M}|(A_{n}\times \Rd) \searrow 0 $ and hence  $| \mu_{j}(A_{n}) | \to 0$. This proves that $\mu_{j}$ is a measure over $\Rd$ ans it is also finite since $C_{M}$ is finite. In addition, from \eqref{Eq:fjDevelop} we have
\begin{equation}
\label{Eq:fjOmuj}
f_{j}(\vec{x}) = \int_{0}^{\vec{x}} \mu_{j}( (-\infty , \vec{u} ] ) d\vec{u} = \mathcal{O}(\mu_{j})(\vec{x}).
\end{equation}
Thus, $f_{j}$ is nothing but $\mathcal{O}(\mu_{j})$, therefore
\begin{equation}
\label{Eq:mujDerfj}
\mu_{j} = \frac{\partial^{2d} f_{j}}{\partial x_{1}^{2} ... \partial x_{d}^{2}}. \quad \blacksquare  
\end{equation}

When $\sigma_{X_{j}}^{2} = 0$ the distribution $\frac{\partial^{2d} f_{j}}{\partial x_{1}^{2} ... \partial x_{d}^{2}}$ is not necessarily a measure, but such case does not really intervene in the decomposition \eqref{Eq:OMexpansionL2nu} ($X_{j} = 0$). For simplicity, we assume from now on that $\sigma_{X_{j}}^{2} > 0$ for all $j \in \mathbb{N}$ (for the case where the sum \eqref{Eq:MvarphiExpansion} is finite we have nothing more to prove).

\begin{lemma}
\label{Lemma:ExpansionMsmooth}
$M$ has the following expansion
\begin{equation}
\label{Eq:MvarphiExpansionSmooth}
\langle M , \varphi \rangle = \sum_{j \in \mathbb{N}}X_{j}\langle \mu_{j} , \varphi \rangle, \quad \forall \varphi \in \DRd.
\end{equation}
\end{lemma}

\textbf{Proof of Lemma \ref{Lemma:ExpansionMsmooth}:} For $\varphi \in \DRd$, one has $\varphi p \in L^{2}(\Rd,\nu)$. Thus, from expansion \eqref{Eq:OMexpansionL2nu} we have

\begin{equation}
\label{Eq:UvarphiL2nu}
\langle \mathcal{O}(M) , \varphi \rangle = \int_{\Rd} \mathcal{O}(M)(x) \varphi(x) p(x) d\nu(x) = \sum_{j \in \mathbb{N}} X_{j} \langle  \nu , f_{j} \varphi p \rangle = \sum_{j \in \mathbb{N}} X_{j} \int_{\Rd}f_{j}(x)\varphi(x)dx.
\end{equation}

Considering the derivative relations \eqref{Eq:OMDerMIntegrals} and \eqref{Eq:mujDerfj}, we conclude 
\begin{equation}
\label{Eq:ExpansionMvarphiSmooth}
\begin{aligned}
\langle M , \varphi \rangle &= \langle \frac{\partial^{2d} \mathcal{O}(M)}{\partial x_{1}^{2} ... \partial x_{d}^{2}} , \varphi \rangle \\
&= \langle \mathcal{O}(M) , \frac{\partial^{2d} \varphi}{\partial x_{1}^{2} ... \partial x_{d}^{2}} \rangle \\
&= \sum_{j \in \mathbb{N}}X_{j}\langle f_{j} ,  \frac{\partial^{2d} \varphi}{\partial x_{1}^{2} ... \partial x_{d}^{2}} \rangle \\
&= \sum_{j \in \mathbb{N}} X_{j} \langle \frac{\partial^{2d} f_{j}}{\partial x_{1}^{2} ... \partial x_{d}^{2}} , \varphi \rangle = \sum_{j \in \mathbb{N}} X_{j} \langle \mu_{j} , \varphi \rangle, \quad \forall \varphi \in \DRd. \quad \blacksquare
\end{aligned}
\end{equation}

The objective of the following lemmas is to extend the expansion \eqref{Eq:MvarphiExpansionSmooth} to $\varphi \in \mathcal{M}_{B}(\Rd)$. We begin with an important covariance to compute. The semi-stochastic Fubini Theorem \ref{Theo:StochasticFubini} will be used here.

\begin{lemma}
\label{Lemma:E(M(A)Xj)}
The following formula holds for every $A \in \BorelRd$:
\begin{equation}
\label{Eq:E(M(A)Xj)}
\mathbb{E}\left( M(A) X_{j} \right) = \sigma_{X_{j}}^{2} \mu_{j}(A).
\end{equation}
\end{lemma}

\textbf{Proof of Lemma \ref{Lemma:E(M(A)Xj)}:} By definition of $\mathcal{O}(M)$ and $X_{j}$ (Eq. \eqref{Eq:ObjectsExpansionOML2nu}) we have
\begin{equation}
\label{Eq:XjReWritten}
\begin{aligned}
X_{j} &= \int_{\Rd} \mathcal{O}(M)(\vec{y})f_{j}(\vec{y})d\nu(\vec{y}) \\
 &= \int_{\Rd} \int_{0}^{\vec{y}} \int_{\Rd} \mathbf{1}_{(-\infty , \vec{u}]}(\vec{s}) dM(\vec{s}) d\vec{u} f_{j}(\vec{y})d\nu(\vec{y}) \\
 &= \int_{\Rd\times\Rd} \int_{\Rd} \mathbf{1}_{(-\infty , \vec{u}]}(\vec{s}) dM(\vec{s}) \theta_{\vec{y}}(\vec{u})f_{j}(\vec{y})d( \ell^{\otimes d}\otimes \nu )(\vec{u} , \vec{y}),
\end{aligned}
\end{equation} 
where $\theta_{\vec{y}}$ is the function $\theta_{\vec{y}} : \Rd \to \lbrace -1 , 0 , 1 \rbrace$ such that $\int_{\Rd} \theta_{\vec{y}}(\vec{u}) \varphi(\vec{u}) d\vec{u} = \int_{0}^{\vec{y}} \varphi(\vec{u}) d\vec{u}$ for every $\varphi \in C(\Rd)$\footnote{The function $\theta_{\vec{y}}$ is just the indicator function of $[0 , \vec{y}]$ when the components of $\vec{y}$ are all positive. When they are not, corresponding minus signs must be added in order to make the integrals coincide. In any case, $\theta_{\vec{y}}$ has compact support.}. We shall apply semi-stochastic Fubini Theorem \ref{Theo:StochasticFubini} to switch integral signs in \eqref{Eq:XjReWritten}. Consider the measure $\lambda$ over $\RdxRd$ given by $d\lambda( \vec{u} , \vec{y} ):=  \theta_{\vec{y}}(\vec{u})f_{j}(\vec{y}) d( \ell^{\otimes d} \otimes \nu )( \vec{u} , \vec{y} )$. $\lambda$ is finite since by (deterministic) Fubini
\begin{equation}
\begin{aligned}
|\lambda|(\RdxRd) = \int_{\RdxRd} | \theta_{\vec{y}}(\vec{u})f_{j}(\vec{y}) |d( \ell^{\otimes d} \otimes \nu )( \vec{u} , \vec{y} ) &= \int_{\Rd} \int_{\Rd} | \theta_{\vec{y}}(\vec{u})| d\vec{u} |f_{j}(\vec{y})| d\nu(\vec{y}) \\
&\leq \int_{\Rd} |y_{1}|... |y_{d}| |f_{j}(\vec{y})| d\nu(\vec{y}) < \infty.
\end{aligned}
\end{equation}
For condition \eqref{It:Fubini1} we use that both $C_{M}$ and $\lambda$ are finite, so
\begin{equation}
\begin{aligned}
\int_{(\RdxRd)\times (\RdxRd) \times (\RdxRd)} \mathbf{1}_{(-\infty , \vec{u}]}(\vec{s})\mathbf{1}_{(-\infty , \vec{v}]}(\vec{t}) &d\left( |C_{M}|\otimes  |\lambda| \otimes |\lambda | \right)\left( (\vec{s},\vec{t}) , (\vec{u},\vec{y}) , (\vec{v},\vec{z}) \right)  \\
&\leq |C_{M}|(\RdxRd)\left[|\lambda|(\RdxRd)\right]^{2} < \infty.
\end{aligned}
\end{equation}
For condition \eqref{It:Fubini2}, we have to study the function
\begin{equation}
\label{Eq:FunctionToUseFubini}
\left(  (\vec{u},\vec{y}) , (\vec{v},\vec{z}) \right) \mapsto \int_{\RdxRd} \mathbf{1}_{(-\infty , \vec{u}]}(\vec{s})\mathbf{1}_{(-\infty , \vec{v}]}(\vec{t}) dC_{M}( \vec{s} , \vec{t} ) = C_{M}\left( (-\infty , \vec{u} ]\times (-\infty , \vec{v}]  \right).
\end{equation}
This function is clearly bounded and it does not depend upon $\vec{y},\vec{z}$, so it is continuous in such components. Moreover, from Proposition \ref{Prop:CMprimitiveContinuous} it also follows that \eqref{Eq:FunctionToUseFubini} is continuous over $ (E\times \Rd)\times(E\times \Rd)$, being $E \in \Borel{\Rd}$ such that $\ell^{\otimes d}(E^{c}) = 0$, and therefore such that $|\lambda|\left( \left[ E \times \Rd \right]^{c} \right) = 0$. Semi-stochastic Fubini Theorem can then be applied to switch the integral order in \eqref{Eq:XjReWritten}, obtaining
\begin{equation}
X_{j} = \int_{\Rd\times\Rd} \int_{\Rd} \mathbf{1}_{(-\infty , \vec{u}]}(\vec{s}) dM(\vec{s}) d\lambda(\vec{u},\vec{y}) = \int_{\Rd} \int_{\Rd\times \Rd} \mathbf{1}_{(-\infty , \vec{u}]}(\vec{s}) d\lambda( \vec{u},\vec{y} ) dM(\vec{s}).
\end{equation}
Using formula \eqref{Eq:CovarianceIntegrals} and (deterministic) Fubini Theorem, we obtain
\begin{equation}
\label{Eq:MAXj}
\begin{aligned}
\mathbb{E}\left(  M(A) X_{j}  \right) &= \mathbb{E}\left( \int_{\Rd} \mathbf{1}_{A}(\vec{t})dM(\vec{t}) \int_{\Rd} \int_{\Rd\times \Rd} \mathbf{1}_{(-\infty , \vec{u}]}(\vec{s}) d\lambda( \vec{u},\vec{y} ) dM(\vec{s}) \right) \\
&= \int_{\RdxRd}  \mathbf{1}_{A}(\vec{t}) \int_{\RdxRd}   \mathbf{1}_{(-\infty , \vec{u}]}(\vec{s}) d\lambda( \vec{u},\vec{y} ) dC_{M}( \vec{t} , \vec{s} ) \\
&=\int_{\RdxRd} \int_{\RdxRd}   \mathbf{1}_{A}(\vec{t})\mathbf{1}_{(-\infty , \vec{u}]}(\vec{s})  dC_{M}( \vec{t} , \vec{s} ) d\lambda( \vec{u},\vec{y} ) \\
&= \int_{\Rd} \int_{0}^{\vec{y}} C_{M}(A \times (-\infty , \vec{u}] ) d\vec{u} f_{j}(\vec{y})d\nu(\vec{y}) \\
&= \int_{\Rd} \mathcal{O}\big( C_{M}(A\times \ \cdot \ ) \big)(\vec{y}) f_{j}(\vec{y})d\nu(\vec{y}) \\
&= \sigma_{X_{j}}^{2}\mu_{j}(A). \quad \blacksquare
\end{aligned}
\end{equation}

\begin{lemma}
\label{Lemma:LambdaWellDef}
The bilinear form $\Lambda : \mathcal{M}_{B}(\Rd)\times \mathcal{M}_{B}(\Rd) \to \mathbb{R}$ given by
\begin{equation}
\label{Eq:DefLambda}
\Lambda(\varphi , \phi) = \sum_{j=1}^{\infty} \sigma_{X_{j}}^{2} \langle \mu_{j} , \varphi \rangle \langle \mu_{j} , \phi \rangle
\end{equation}
is well-defined, the series being absolutely convergent.
\end{lemma}

\textbf{Proof of Lemma \ref{Lemma:LambdaWellDef}:} Let us define the sequence of finite random measures
\begin{equation}
\label{Eq:DefMn}
M_{n} = \sum_{j \leq n} X_{j} \mu_{j}, \quad n \in \mathbb{N}.
\end{equation}
Their covariance measures $C_{M_{n}}$ are given by
\begin{equation}
\label{Eq:CMn}
\begin{aligned}
C_{M_{n}}(A\times B) &= \mathbb{E} ( \ \sum_{j\leq n} \sum_{k \leq n} X_{j} X_{k} \mu_{j}(A)\mu_{k}(B) \ ) \\
&= \sum_{j \leq n} \sigma_{X_{j}}^{2} \mu_{j}(A) \mu_{j}(B),
\end{aligned}
\end{equation}
where we have used $\mathbb{E}(X_{j}X_{k} ) = \sigma_{X_{j}}^{2}\delta_{j,k}$. In addition, using Lemma \ref{Lemma:E(M(A)Xj)} we conclude
\begin{equation}
\label{Eq:CMnMdevelop}
\begin{aligned}
\mathbb{E}\left( M(A)M_{n}(B) \right) &= \mathbb{E}( \ \sum_{j \leq n} M(A)X_{j}\mu_{j}(B)  \ ) \\
&= \sum_{j \leq n} \mathbb{E}\left( M(A)X_{j}  \right) \mu_{j}(B) \\
&= \sum_{j \leq n} \sigma_{X_{j}}^{2}\mu_{j}(A)\mu_{j}(B) &= C_{M_{n}}(A\times B).
\end{aligned}
\end{equation}
Developing the expression $\mathbb{E}\left( \big( M(A)-M_{n}(A) \big) \big( M(B) - M_{n}(B) \big) \right)$, one concludes from \eqref{Eq:CMn} and \eqref{Eq:CMnMdevelop} that $M-M_{n}$ is an m-cov  finite random measure with covariance
\begin{equation}
\label{Eq:CMminusMn}
C_{M-M_{n}} = C_{M} - C_{M_{n}} = C_{M} - \sum_{j \leq n} \sigma_{X_{j}}^{2}\mu_{j}\otimes \mu_{j}.
\end{equation}
Since $C_{M-M_{n}}$ is a finite covariance measure it must be positive definite, having
\begin{equation}
\langle C_{M}  - C_{M_{n}} , \varphi \otimes \varphi \rangle =  \langle C_{M-M_{n}} , \varphi \otimes \varphi \rangle \geq 0, \quad \forall \varphi \in \mathcal{M}_{B}(\Rd),
\end{equation}
which implies
\begin{equation}
\label{Eq:InequalityCMnCM}
\langle C_{M_{n}} , \varphi \otimes \varphi \rangle \leq \langle C_{M} , \varphi \otimes \varphi \rangle, \quad \forall \varphi \in \mathcal{M}_{B}(\Rd).
\end{equation}
Using Cauchy-Schwartz inequality, we conclude for every $n$
\begin{equation}
\label{Eq:CMnConvergentSeries}
\begin{aligned}
\sum_{j \leq n}  \sigma_{X_{j}}^{2}|\langle \mu_{j} , \varphi \rangle | |\langle \mu_{j} , \phi \rangle | &\leq \sqrt{\sum_{j \leq n} \sigma_{X_{j}}^{2} | \langle \mu_{j} , \varphi \rangle |^{2} }\sqrt{\sum_{j \leq n} \sigma_{X_{j}}^{2} | \langle \mu_{j} , \phi \rangle |^{2} } \\
&= \sqrt{\langle C_{M_{n}} , \varphi \otimes \varphi \rangle }\sqrt{\langle C_{M_{n}} , \phi \otimes \phi \rangle } \\
&\leq \sqrt{\langle C_{M} , \varphi \otimes \varphi \rangle }\sqrt{\langle C_{M} , \phi \otimes \phi \rangle } < \infty,
\end{aligned}
\end{equation}
which proves that the series \eqref{Eq:DefLambda} is absolutely convergent and thus the bilinear form $\Lambda$ is well-defined. $\blacksquare$

The following Lemma is the crucial part where an argument essentially different as those found in the proof of the classical KL expansion is needed. here, an adequate use of Lusin's Theorem will help us to conclude the convergence against measurable and bounded functions.

\begin{lemma}
\label{Lemma:LambdaCMvarphi}
The following equality holds
\begin{equation}
\label{Eq:LambdaCMvarphi}
\Lambda(\varphi , \varphi) = \langle C_{M} , \varphi \otimes \varphi \rangle, \quad \forall \varphi \in \mathcal{M}_{B}(\Rd).
\end{equation}
\end{lemma}

\textbf{Proof of Lemma \ref{Lemma:LambdaCMvarphi}:} 
Lemma \ref{Lemma:ExpansionMsmooth} guarantees that \eqref{Eq:LambdaCMvarphi} holds for $\varphi \in \DRd$. We will extend it to $\varphi \in \mathcal{M}_{B}(\Rd)$. We begin by considering $\varphi \neq 0 $ of the form $\varphi = \mathbf{1}_{I}$, where $I \subset \R^{d}$ is a rectangle $ I = I_{1} \times ... \times I_{d}$, each $I_{j}$ being an interval of $\R$. In such case $\varphi$ can be approximated point-wisely by a sequence of functions in $\DRd$, the sequence being dominated by $\| \varphi \|_{\infty}$. Let $\epsilon > 0 $. Since $C_{M}$ is a finite measure, from dominated convergence we can choose $\phi \in \DRd$ approaching $\varphi$ so that $\| \phi \|_{\infty} \leq \| \varphi \|_{\infty}$  and so that
\begin{equation}
\label{Eq:BoundedCMVarphiDiffPhi}
\langle |C_{M}| , |\varphi - \phi |\otimes |\varphi - \phi | \rangle < \frac{\epsilon^{2}}{64\|\varphi \|_{\infty}^{2}|C_{M}|(\RdxRd)}.
\end{equation}
Now, from triangular inequality we have
\begin{equation}
\label{Eq:BoundingLambdaCM}
\begin{aligned}
\left| \Lambda(\varphi , \varphi ) - \langle C_{M} , \varphi \otimes \varphi \rangle \right| \leq &\underbrace{| \  \Lambda(\varphi , \varphi ) -  \sum_{j \leq n} \sigma_{X_{j}}^{2}| \langle \mu_{j} , \varphi \rangle|^{2} |  \ }_{(a)} + \underbrace{| \    \sum_{j \leq n} \sigma_{X_{j}}^{2}| \langle \mu_{j} , \varphi \rangle|^{2} -  \sum_{j \leq n} \sigma_{X_{j}}^{2}| \langle \mu_{j} , \phi \rangle|^{2}  | \ }_{(b)} \\
&+ \underbrace{ \ |   \sum_{j \leq n} \sigma_{X_{j}}^{2}| \langle \mu_{j} , \phi \rangle|^{2} - \langle C_{M} , \phi \otimes \phi \rangle | \ }_{(c)}+ \underbrace{ | \  \langle C_{M} , \phi \otimes \phi \rangle   - \langle C_{M} , \varphi\otimes \varphi \rangle  | \ }_{(d)}.
\end{aligned}
\end{equation}
By Lemmas \ref{Lemma:LambdaWellDef} and \ref{Lemma:ExpansionMsmooth}, there exists $n_{0}$ such that both terms $(a)$ and $(c)$ are smaller than $\frac{\epsilon}{4}$ if $n \geq n_{0}$. For the term $(d)$ we use the symmetry of $C_{M}$, the Cauchy-Schwarz inequality and inequality \eqref{Eq:BoundedCMVarphiDiffPhi} to obtain
\begin{equation}
\label{Eq:Bounding(d)}
\begin{aligned}
(d) &= | \langle C_{M} , \varphi\otimes \varphi  - \phi \otimes \phi \rangle | \\
&=  | \langle  C_{M}  , (\varphi + \phi)  \otimes ( \varphi  - \phi ) \rangle | \\
&\leq \sqrt{\langle C_{M} , (\varphi + \phi)\otimes (\varphi + \phi) \rangle } \sqrt{\langle C_{M} , (\varphi - \phi)\otimes (\varphi - \phi) \rangle } \\
&< \sqrt{ |C_{M}|(\RdxRd) 4 \| \varphi \|_{\infty}^{2}  } \sqrt{  \frac{\epsilon^{2}}{64\|\varphi \|_{\infty}^{2}|C_{M}|(\RdxRd)}  } \quad  = \frac{\epsilon}{4}.
\end{aligned}
\end{equation}
On the other hand, for the term $(b)$ we can do similarly, considering the covariance measure $C_{M_{n}}$ (Eq. \eqref{Eq:CMn}) and inequality \eqref{Eq:InequalityCMnCM}:
\begin{equation}
\label{Eq:Bounding(b)}
\begin{aligned}
(b) &= | \langle C_{M_{n}} , \varphi\otimes \varphi  - \phi \otimes \phi \rangle | \\
&=  | \langle  C_{M_{n}}  , (\varphi + \phi)  \otimes ( \varphi  - \phi ) \rangle | \\
&\leq \sqrt{\langle C_{M_{n}} , (\varphi + \phi)\otimes (\varphi + \phi) \rangle } \sqrt{\langle C_{M_{n}} , (\varphi - \phi)\otimes (\varphi - \phi) \rangle } \\
&\leq \sqrt{\langle C_{M} , (\varphi + \phi)\otimes (\varphi + \phi) \rangle } \sqrt{\langle C_{M} , (\varphi - \phi)\otimes (\varphi - \phi) \rangle }  < \frac{\epsilon}{4}.  \\
\end{aligned}
\end{equation}
We conclude that $\left| \Lambda(\varphi , \varphi) - \langle C_{M} , \varphi \otimes \varphi \rangle  \right| \leq  \epsilon$ for every $\epsilon > 0 $, and therefore $
\Lambda(\varphi , \varphi ) = \langle C_{M} , \varphi \otimes \varphi \rangle$ for every $\varphi$ of the form $\varphi = \mathbf{1}_{I}$. By bi-linearity of $\Lambda$, we can easily extend this result to every $\varphi$ in the space
\begin{equation}
\label{Eq:DefSpaceE}
\mathcal{E} := \Span \lbrace  \ \mathbf{1}_{I} \ \big| \ I \subset \Rd \hbox{ rectangle }  \rbrace.
\end{equation}
Now, in order to extend this result to any $\varphi \in \mathcal{M}_{B}(\Rd)$, we use Lusin's Theorem \ref{Theo:Lusin} applied to the space $\Rd$ with the finite measure $|C_{M}|( \ \cdot \times \Rd)$. Given $\varphi \in \mathcal{M}_{B}(\Rd)$, $\varphi \neq 0 $, and given $\epsilon > 0 $, there exists a closed set $E \subset \Rd$ such that $\varphi$ is continuous over $E$ (with the subspace topology) and such that
\begin{equation}
\label{Eq:BoundCMEcLusin}
|C_{M}|(E^{c} \times \Rd) < \frac{\epsilon^{2}}{ 1536\| \varphi \|^{4}_{\infty}|C_{M}|(\RdxRd)  }.
\end{equation}
Consider a typical Riemann-alike approximation of $\varphi$, done through a sequence of functions of the form
\begin{equation}
\label{Eq:DefPhin}
\phi_{n} := \sum_{j = 1}^{n}\varphi( x_{j}^{n} ) \mathbf{1}_{I_{j}^{n}},
\end{equation}
where for each $n$, $(I_{j}^{n})_{j = 1 , ... , n} $ is a collection of rectangles forming a partition of a subset $K_{n}$ of $\Rd$, satisfying that $K_{n} \nearrow \Rd$ and $\displaystyle\max_{j=1,...,n}\diam(I_{j}^{n}) \to 0 $ as $n \to \infty$; and $x_{j}^{n} \in I_{j}^{n}$ is a tag-point chosen so $x_{j}^{n} \in E$ when $I_{j}^{n} \cap E \neq \emptyset$. Let $x \in E$. Denote $j_{n}$ the index of the interval $I_{j}^{n}$ where $x$ belongs to. By construction the sequence $(x_{j_{n}}^{n})_{n}$ is in $E$, and $x_{j_{n}}^{n} \to x \in E$ as $n \to \infty$. Since $\varphi$ is continuous over $E$ with the subspace topology, we have $\phi_{n}(x) = \varphi(x_{j_{n}}^{n}) \to \varphi(x)$ as $n \to \infty$. In addition, one has $\| \phi_{n} \|_{\infty} \leq \| \varphi \|_{\infty}$. The sequence \eqref{Eq:DefPhin} converges thus point-wisely and dominated to $\varphi$ over $E$. By dominated convergence ($C_{M}$ is a measure), for every $\epsilon >0 $ we can select $\phi$ of the form \eqref{Eq:DefPhin} (so in $\mathcal{E}$) with $\| \phi \|_{\infty} \leq \| \varphi \|_{\infty}$ such that
\begin{equation}
\int_{E\times E} |\varphi - \phi | \otimes | \varphi - \phi | d|C_{M}| < \frac{\epsilon^{2}}{128\| \varphi \|_{\infty}^{2} |C_{M}|(\RdxRd)  }.
\end{equation}
By splitting integrals, using elementary bounds and the symmetry of $|C_{M}|$, we have
\begin{equation}
\label{Eq:BoundingDoubleIntegralOverExE}
\begin{aligned}
\langle |C_{M}| , |\varphi - \phi|\otimes |\varphi - \phi| \rangle &= \int_{E\times E}|\varphi - \phi | \otimes | \varphi - \phi | d|C_{M}| + \int_{E\times E^{c}}|\varphi - \phi | \otimes | \varphi - \phi | d|C_{M}|  \\
&\ \ \ + \int_{E^{c}\times E}|\varphi - \phi | \otimes | \varphi - \phi | d|C_{M}| + \int_{E^{c}\times E^{c}}|\varphi - \phi | \otimes | \varphi - \phi | d|C_{M}| \\
&\leq  \frac{\epsilon^{2}}{128\| \varphi \|_{\infty}^{2}|C_{M}|(\RdxRd) }  +  3 \cdot 4 \| \varphi \|^{2}_{\infty} |C_{M}|(E^{c}\times \Rd)    \\
&<   \frac{\epsilon^{2}}{128\| \varphi \|_{\infty}^{2}|C_{M}|(\RdxRd) }  +  12 \| \varphi \|^{2}_{\infty}  \frac{\epsilon^{2}}{ 1536\| \varphi \|^{4}_{\infty}|C_{M}|(\RdxRd)  } \\
&= \frac{\epsilon^{2}}{64 \| \varphi \|_{\infty}^{2} |C_{M}|(\RdxRd) }. 
\end{aligned}
\end{equation}
With this set up, we can study the expression $| \Lambda(\varphi , \varphi ) - \langle C_{M} , \varphi \otimes \varphi \rangle |$ for any $\varphi \in \mathcal{M}_{B}(\Rd)$ by using \textit{the same} splitting arguments exposed in \eqref{Eq:BoundingLambdaCM}, using $\phi \in \mathcal{E}$ constructed as above. Expressions $(a)$ and $(c)$ can be bounded by $\frac{\epsilon}{4}$ for $n \geq n_{0}$ for some $n_{0}$. In expressions $(b)$ and $(d)$ we can also follow line by line the arguments \eqref{Eq:Bounding(d)} and \eqref{Eq:Bounding(b)} to bound both of them by $\frac{\epsilon}{4}$. We conclude once again that $\left| \Lambda(\varphi , \varphi) - \langle C_{M} , \varphi \otimes \varphi \rangle  \right| \leq  \epsilon$ for every $\epsilon > 0 $, and therefore, since $\varphi$ is arbitrary, 
\begin{equation}
 \Lambda(\varphi , \varphi ) = \langle C_{M} , \varphi \otimes \varphi \rangle = \sum_{j \in \mathbb{N}} \sigma_{X_{j}}^{2}|\langle \mu_{j} , \varphi \rangle|^{2}, \quad \forall \varphi \in \mathcal{M}_{B}(\Rd). \quad \blacksquare
\end{equation}

Now we finish the proof of Theorem \ref{Theo:KLMeasure}.

\textbf{Proof of Theorem \ref{Theo:KLMeasure}:} Everything being set up, we use Lemmas \ref{Lemma:LambdaWellDef}, \ref{Lemma:LambdaCMvarphi} and expression \eqref{Eq:CMminusMn} to obtain
\begin{equation}
\begin{aligned}
\mathbb{E}\big( | \langle M , \varphi \rangle - \sum_{j \leq n} X_{j}\langle \mu_{j} , \varphi \rangle |^{2} \big) &= \mathbb{E}( | \langle M  - M_{n}, \varphi \rangle |^{2} ) \\
&= \langle C_{M-M_{n}} , \varphi \otimes \varphi \rangle \\
&= \langle C_{M} - C_{M_{n}} , \varphi \otimes \varphi \rangle \\
&= \Lambda( \varphi , \varphi ) - \sum_{j \leq n} \sigma_{X_{j}}^{2} |\langle \mu_{j} , \varphi \rangle |^{2} \to 0, \quad \forall \varphi \in \mathcal{M}_{B}(\Rd).\quad \blacksquare
\end{aligned}
\end{equation}

From Theorem \ref{Theo:KLMeasure} it follows that the covariance measure $C_{M}$ has the expansion
\begin{equation}
\label{Eq:ExpansionCM}
\langle C_{M} , \varphi \otimes \phi \rangle = \sum_{j \in \mathbb{N}} \sigma_{X_{j}}^{2}\langle \mu_{j} , \varphi \rangle \langle \mu_{j} , \phi \rangle, \quad \forall \varphi,\phi \in \mathcal{M}_{B}(\Rd).
\end{equation}
There is a slightly stronger convergence mode for this expansion: if we fix $\varphi$ (or $\phi$) then the measure $\langle C_{M_{n}} , \varphi \otimes \ \cdot \ \rangle $ converges in absolute variation to $\langle C_{M} , \varphi \otimes \ \cdot \ \rangle$.
\begin{proposition}
\label{Prop:ExpansionCMabsVar}
For every $\varphi \in \mathcal{M}_{B}(\Rd)$, one has
\begin{equation}
\big|  \langle C_{M} , \varphi \otimes \ \cdot \ \rangle - \sum_{j \leq n} \sigma_{X_{j}}^{2}\langle \mu_{j} , \varphi \rangle \mu_{j} \  \big|(\Rd) \to 0, \quad \hbox{ as } n \to \infty.
\end{equation}
\end{proposition}

\textbf{Proof of Proposition \ref{Prop:ExpansionCMabsVar}:} Let $\varphi \in \mathcal{M}_{B}(\Rd)$. We use expression \eqref{Eq:TotalVariationAlternative} for the total-variation measure of $\langle C_{M-M_{n}} , \varphi \otimes \ \cdot \  \rangle $. Since \eqref{Eq:CMminusMn} implies $\langle C_{M-M_{n}} , \phi \otimes \phi \rangle \leq \langle C_{M} , \phi\otimes\phi \rangle$ for every $\phi \in \mathcal{M}_{B}(\Rd)$,  we have
\begin{equation}
\begin{aligned}
\left| \langle C_{M}-C_{M_{n}} , \varphi \otimes \ \cdot \ \rangle \right|(\Rd) &=\left| \langle C_{M-M_{n}} , \varphi \otimes \ \cdot \ \rangle \right|(\Rd) \\ 
&= \sup_{\phi \in \mathcal{M}(\Rd), |\phi | = 1} \left| \langle  C_{M-M_{n}} , \varphi \otimes \phi \rangle \right| \\
&\leq \sup_{\phi \in \mathcal{M}(\Rd), |\phi | = 1}  \sqrt{ \langle C_{M-M_{n}} , \varphi \otimes \varphi \rangle} \sqrt{ \langle C_{M-M_{n}} , \phi \otimes \phi \rangle }  \\
&\leq \sqrt{ \langle C_{M-M_{n}} , \varphi \otimes \varphi \rangle} \sup_{\phi \in \mathcal{M}(\Rd), |\phi | = 1}  \sqrt{ \langle C_{M} , \phi \otimes \phi \rangle }  \\
&\leq \sqrt{ \langle C_{M-M_{n}} , \varphi \otimes \varphi \rangle} \sqrt{ |C_{M}|(\RdxRd) }  \xrightarrow[n \to \infty]{} 0. \quad \blacksquare
\end{aligned}
\end{equation}

To finish, we present the expansion of function-regulated random measures, which covers some non-finite random measures cases (such as White Noise or Poisson processes with non-finite intensity). Let us introduce the following definition.

\begin{definition}
\label{Def:MRegulatedf}
A measure $\mu \in \mathscr{M}(\Rd)$ (resp., an m-cov random measure $M$ over $\Rd$) is said to be function-regulated if there exists a strictly positive and locally bounded measurable function $f$ such that $\mu/f$ is finite (resp., such that $M/f$ is finite).
\end{definition}

It is not difficult to verify that $M$ is regulated by $f$ if and only if $C_{M}$ is regulated by $f \otimes f$.

\begin{theorem}
\label{Theo:KLExpansionMnonfinite}
Let $M$ be an m-cov random measure over $\Rd$ regulated by a function $f$. Then, there exist a sequence of pairwise uncorrelated random variables with summable variances $(X_{n})_{n \in \mathbb{N}}$ and a linearly independent sequence of real measures $(\mu_{n})_{n \in \mathbb{N}} \subset \mathscr{M}(\Rd)$, all of them regulated by $f$, such that
\begin{equation}
\label{Eq:MvarphiExpansionNonFinite}
\langle M , \varphi \rangle = \sum_{n \in \mathbb{N}}X_{n}\langle \mu_{n} , \varphi \rangle, \quad \forall \varphi \in \mathcal{M}_{B,c}(\Rd),
\end{equation}
the series being considered in a mean-square sense. In addition, the covariance measure $C_{M}$ satisfies
\begin{equation}
\label{Eq:CMExpansionConvergenceNonFinite}
| \ \langle C_{M} , \varphi \otimes \ \cdot \ \rangle -  \sum_{j\leq n} \sigma_{X_{j}}^{2}\langle \mu_{j} , \varphi \rangle  \mu_{j} \ |(K)  \xrightarrow[n \to \infty]{} 0, \quad \forall \varphi \in \mathcal{M}_{B,c}(\Rd), \forall K \subset \Rd \hbox{ compact}.
\end{equation}
\end{theorem}

\textbf{Proof of Theorem \ref{Theo:KLExpansionMnonfinite}:} We apply Theorem \ref{Theo:KLMeasure} to the finite random measure $\frac{1}{f}M$, obtaining thus
\begin{equation}
\langle \frac{1}{f}M , \phi \rangle = \sum_{j \in \mathbb{N}} X_{j} \langle \nu_{j} , \phi \rangle, \quad \phi \in \mathcal{M}_{B}(\Rd),
\end{equation}
with $\nu_{j} \in \mathscr{M}_{F}(\Rd)$ for every $j$. By posing $\mu_{j} = f\nu_{j}$ and using $\langle M , \varphi \rangle = \langle \frac{1}{f}M , f\varphi \rangle$ the result follows. The convergence for the covariance \eqref{Eq:CMExpansionConvergenceNonFinite} is obtained following the same arguments as in Proposition \ref{Prop:ExpansionCMabsVar}. $\blacksquare$

\section{Concluding remarks}
\label{Sec:Remarks}

The following remarks are meant to clarify some important points, expose some remarkable cases and to provide ideas of extensions and applications of the results here obtained.

\subsection{Hilbert space-based approach}
\label{Sec:HilbertApproach}

As a general rule, KL expansions are constructed from a Hilbert space perspective as presented in Section \ref{Sec:KLHilbert}. Our case is no exception, although the Hilbert space in game is not yet explicitly shown. To make it precise, we consider the space of finite measures $\mathscr{M}_{F}(\Rd)$ endowed with the following bilinear positive-definite form
\begin{equation}
\label{Eq:DefineHilberProdForMeasures}
(\mu_{1} , \mu_{2} )_{E} \mapsto \left( \mathcal{O}(\mu_{1}) , \mathcal{O}( \mu_{2} ) \right)_{L^{2}(\Rd,\nu)},
\end{equation}
where $\nu \in \mathscr{M}_{F}(\Rd)$ is defined as in \eqref{Eq:Defnu}. By construction,  $\mathcal{O} : \mathscr{M}_{F}(\Rd) \to L^{2}(\Rd , \nu) $. From the continuity of $\mathcal{O}(\mu)$ and since $\nu$ has a density, one has $\int_{\Rd}|\mathcal{O}(\mu)|^{2}d\nu  = 0 \Longleftrightarrow \mathcal{O}(\mu) = 0 $. Since $\mathcal{O}(\mu)$ is the primitive of $\vec{u} \mapsto \mu( (-\infty , \vec{u}] )$,  $\mathcal{O}(\mu) = 0$ implies $\mu( (-\infty , \vec{u}] ) = 0$ $\vec{u}$-almost-everywhere, and from right-continuity it must be null. We conclude $(\mu , \mu )_{E} = 0 \Longleftrightarrow \mu = 0$, and thus \eqref{Eq:DefineHilberProdForMeasures} is a Hilbert product. The completion of $\mathscr{M}_{F}(\Rd)$ with this product is then an abstract separable Hilbert space $E$ in which $\mathscr{M}_{F}(\Rd)$ is dense. This space consists of distributions which are derivatives of $2d$ order of elements in $L^{2}( \Rd , \nu )$. The KL expansion \eqref{Eq:MvarphiExpansion} can be obtained as the Hilbert-based expansion of $M$ with respect to $E$. The measures $(\mu_{n})_{n}$ are orthonormal in $E$, since $f_{n} = \mathcal{O}(\mu_{n})$ (Eq. \eqref{Eq:fjOmuj}). Following Section \ref{Sec:KLHilbert}, one can identify $M$ as a process linearly indexed on $E$ by
\begin{equation}
M(\mu) = \int_{\Rd} \mathcal{O}(M) \mathcal{O}(\mu) d\nu, \quad \forall \mu \in E.
\end{equation}

It is important to remark that \textit{a} KL expansion is always done with respect to a Hilbert space which is chosen with some arbitrariness. For example, in the classical case of a mean-square continuous process over a compact interval one could use $L^{2}([a,b] , \alpha  )$ as basis Hilbert space with some measure with density $\alpha$ rather than $L^{2}([a,b])$, obtaining a strictly different expansion. Similarly, if the process $X$ has a twice-differentiable covariance function, one could use the Sobolev space $H^{1}( (a,b) )$ as reference Hilbert space, obtaining in general a different decomposition than using $L^{2}([a,b])$. Independently of the Hilbert space used, the real contribution of this work is the proof of a stronger convergence than in the Hilbert space sense, here translated in a $\mathcal{M}_{B,c}(\Rd)^{*}$-weak-sense implying set-wise convergence, analogously to the classical uniform-mean-square convergence for the expansion of a mean-square continuous process over a compact interval.

\subsection{Uniqueness of the expansion}
\label{Sec:Uniqueness}

In many senses, the expansion \eqref{Eq:MvarphiExpansion} for a given $M$ is far from being unique. The first source of arbitrariness comes from the choice of $\nu$ in \eqref{Eq:Defnu}, where other enough-fast-decreasing density measures could be chosen. Another arbitrary choice is the use of the anti-derivative operator $\mathcal{O}$ to link $\mathscr{M}_{F}(\Rd)$ to a Hilbert space. Other regularising operators could be used, for instance operators of the form $(1 - \Delta)^{-\alpha}$ with $\alpha > 0$ big enough and $\Delta$ the Laplacian, as it is commonly used for the definition of negative Sobolev spaces. As explained in Section \ref{Sec:HilbertApproach}, the choice of the basis Hilbert space is also arbitrary.

Therefore, the actual question that has here been answered positively for m-cov function-regulated random measures is: \textit{is it possible for a random measure to be expressed as series of deterministic measures weighted by uncorrelated random variables, with an adequate convergence with respect to its measure structure?} The uniqueness of the series is not studied here. One interesting question that arises is if all these constructions have a common reference property which would allow us to speak about ``\textit{the}'' KL expansion of a random measure $M$. For instance, it is expected that if we change $\nu$ for another measure with strictly positive density, then the newly obtained measures $(\mu_{j})_{j \in \mathbb{N}}$ in the KL expansion will be absolutely continuous with respect to the ones obtained with $\nu$.

\subsection{Gaussian case}
\label{Sec:MGauss}

If $M$ is Gaussian, that is, if $(M(A_{1}) , \ldots , M(A_{n}) )$ is a Gaussian vector for every $A_{1},\ldots , A_{n} \in \BoundedBorel{\Rd}$, then the variables $X_{j}$ in the expansion of $M$ are independent and Gaussian, since all the variables involved are constructed linearly. In addition, the convergence of the series also holds almost-surely. This follows from a classical result on almost surely convergence of series of independent random variables with variances forming a convergent series \citep[Section 12.2]{williams1990probability}. In our case,
\begin{equation}
\sum_{n \in \mathbb{N}} \Var( X_{j}\langle \mu_{n} , \varphi \rangle ) = \sum_{n \in \mathbb{N}} \sigma_{X_{n}}^{2} |\langle \mu_{n} , \varphi \rangle|^{2} = \langle C_{M} , \varphi \otimes \varphi \rangle < \infty \ \Longrightarrow \ \sum_{n \in \mathbb{N}} X_{n} \langle \mu_{n} , \varphi \rangle \hbox{ converges a.s.} 
\end{equation}

\subsection{Expansions for non-regular processes and trawl processes.}

Theorems \ref{Theo:KLMeasure} and \ref{Theo:KLExpansionMnonfinite} can be used to obtain, as corollaries, diverse forms of KL expansions of non-regular stochastic processes over $\Rd$. Consider for instance a process of the form $Z(\vec{x}) =  \langle M , \phi_{\vec{x}} \rangle$, with $M$ being a function-regulated m-cov random measure and $\phi_{\vec{x}} \in \mathcal{M}(\Rd)$ for every $\vec{x}$, such that any of the Theorems \ref{Theo:KLMeasure} or \ref{Theo:KLExpansionMnonfinite} holds for any $\vec{x}$. Then, the following traceable KL expansion for $Z$ holds:
\begin{equation}
\label{Eq:ExpansionZRemark}
Z(\vec{x}) = \sum_{n \in \mathbb{N}} X_{n} g_{n}(\vec{x}),
\end{equation} 
where $g_{n}(\vec{x}) = \langle \mu_{n}, \phi_{\vec{x}} \rangle $. Note that the convergence \eqref{Eq:ExpansionZRemark} is mean-square-point-wise, and this holds without requiring any particular regularity on $Z$ (no mean-square continuity or measurability), since $\vec{x}$ is acting just as an index parameter for $\phi_{\vec{x}}$ with respect to which no regularity is required. 

One important example of so-defined stochastic processes are trawl processes \citep{veraart2019modeling,sauri2022nonparametric}. As mentioned in Section \ref{Sec:WhiteNoise}, a Lévy basis with finite variance is an example of m-cov random measure. If $M$ is such a Lévy basis over $[0,\infty)$, then a process defined as
\begin{equation}
Z(t) = \langle M , \mathbf{1}_{A_{t}} \rangle, \quad t \geq 0,
\end{equation}
is called a \textit{trawl} process. For every $t$, $A_{t} \in \Borel{[0,\infty)}$ is called the trawl set. If $M$ is finite, or if $M$ is function-regulated and $A_{t}$ is bounded for every $t$, then expansion \eqref{Eq:ExpansionZRemark} holds, providing thus a KL expansion for trawl processes.

\subsection{General m-cov random measure case}
\label{Sec:RemGeneralMcov}

It is not clear if a general m-cov random measure $M$ can be regulated by a function $f$ as it is required in Theorem \ref{Theo:KLExpansionMnonfinite}. One thing that can always be done is to construct a KL decomposition \textit{locally}. That is, for every $D \in \BoundedBorel{\Rd}$, $M$ has a decomposition of the form \eqref{Eq:MvarphiExpansion} for every $\varphi \in \mathcal{M}_{B,c}(\Rd)$ null outside $D$. This holds since one can focus on the compactly supported random measure $\mathbf{1}_{D}M$, which is finite and thus Theorem \ref{Theo:KLMeasure} applies. In such case, the measures $\mu_{j}$ and the random variables $X_{j}$ depend upon the set $D$.

\begin{appendix}

\section{Proofs of claims presented in Section \ref{Sec:RandomMeasures}}
\label{App:ProofsSectionRM}

\subsection{Proof of Proposition \ref{Prop:MsigmaAdditive}}

Let $M$ be an m-cov random measure over $\Rd$ with covariance measure $C_{M}$. Let $(A_{n})_{n \in \mathbb{N}} $ be a sequence of pairwise disjoint bounded Borel subsets of $\Rd$ such that $\bigcup_{n \in \mathbb{N}}A_{n} \in \BoundedBorel{\Rd}$. By symmetry of $C_{M}$, we have
\begin{equation}
\label{Eq:MsigmaAdditiveProof}
\mathbb{E}\Big( | M\big( \bigcup_{j \in \mathbb{N}}A_{j} \big) - \sum_{j \leq n} M(A_{j})  |^{2} \Big) =  C_{M}\Big(  \bigcup_{j \in \mathbb{N}}A_{j}  \times  \bigcup_{j \in \mathbb{N}}A_{j}  \Big) - 2 C_{M}\Big( \bigcup_{j \leq n}A_{j} \times  \bigcup_{j \in \mathbb{N}}A_{j}  \Big) + C_{M}\Big(   \bigcup_{j \leq n}A_{j}  \times  \bigcup_{j \leq n}A_{j}   \Big).  
\end{equation}
Since $C_{M}$ is a measure, by $\sigma$-additivity \eqref{Eq:MsigmaAdditiveProof} must go to $0$ as $n \to \infty$. $\blacksquare$

\subsection{Proof of semi-stochastic Fubini Theorem \ref{Theo:StochasticFubini}}
\label{App:FubiniProof}

The first issue with semi-stochastic Fubini Theorem \ref{Theo:StochasticFubini} is the proper definition of the iterated integrals in \eqref{Eq:FubiniTheorem}. Namely, we require a canonical definitions of a stochastic process which is almost-everywhere mean-square continuous with respect to a deterministic measure (right side of \eqref{Eq:FubiniTheorem}). Here we follow an approach using classical Riemann sums, as it is exposed for example in [Section 4.5]\citep{soong1973random} for $d=1$ over a compact interval (there the measure is used through its bounded variation primitive). Here we need slightly more generality, so we develop explicitly such general definition, but the procedure is essentially the same as in \citep{soong1973random}. We remark that a standard method for defining integrals of stochastic processes with respect to deterministic measures is using the Bochner integral \citep{diestel1974vector}. However, Bochner integrability requires actually stronger conditions than the one required here, therefore we do not follow such approach.\footnote{For information, for defining a stochastic integral of the form $\int_{\Rd}Z(x)d\mu(x)$ with $\mu \in \mathscr{M}(\Rd)$ and $Z$ a second-order process in the spirit of the Bochner integral, one requires, at least, 
\begin{equation}
\label{Eq:ConditionBochner}
\int_{\Rd}\sqrt{C_{Z}(x,x)}d|\mu|(x) < \infty,
\end{equation}
where $C_{Z}$ is the covariance function of $Z$. Such a condition is actually stronger than the one we require, namely (see Lemma \ref{Lemma:IntZmuUnbounded})
\begin{equation}
\label{Eq:ConditionRiemann}
\int_{\RdxRd} |C_{Z}|d|\mu|\otimes|\mu| < \infty.
\end{equation}
That \eqref{Eq:ConditionBochner} implies \eqref{Eq:ConditionRiemann} follows from the Cauchy-Schwarz inequality. An example in which \eqref{Eq:ConditionRiemann} holds but \eqref{Eq:ConditionBochner} does not, is $d=1$, $\mu$ the Lebesgue measure, and
\begin{equation}
C_{Z}(x,y) = \frac{e^{-|x-y|}}{xy} \mathbf{1}_{[1,\infty)}(x)\mathbf{1}_{[1,\infty)}(y).
\end{equation} }

\begin{lemma}[\textbf{Dominated convergence for double sequences}]
\label{Lemma:LDCTdoubleSequence}
Let $(E,\mathscr{E}, \mu)$ be a measure space $(\mu \geq 0)$. Let $(f_{n,m})_{n,m \in \mathbb{N}}$ be a double-sequence of complex functions such that $\sup_{n,m \in \mathbb{N}} |f_{n,m}| \in \mathscr{L}^{1}(E,\mathscr{E},\mu)$. Suppose that the double limit $\lim_{n,m \to \infty} f_{n,m} $ exists $\mu$-almost everywhere. Then, the $\mu-$almost everywhere defined function $f(x) = \lim_{n,m \to \infty} f_{n,m}(x)$ is in $\mathscr{L}^{1}(E,\mathscr{E},\mu)$ and
\begin{equation}
\label{Eq:L1Convergefnmf}
\lim_{n,m \to \infty}\int_{E}|f_{n,m}-f|d\mu = 0.
\end{equation}
\end{lemma}

\textbf{Proof of Lemma \ref{Lemma:LDCTdoubleSequence}:}\footnote{The essential of this proof has been borrowed from the StackExchange discussion \url{https://math.stackexchange.com/questions/448931/dominated-convergence-thm-dct-for-double-sequences}, consulted for the last time March the 27th 2025.} $f \in \mathscr{L}^{1}(E,\mathscr{E},\mu)$ since $\sup_{n,m \in \mathbb{N}} |f_{n,m}| \in \mathscr{L}^{1}(E,\mathscr{E},\mu)$. Suppose the double limit \eqref{Eq:L1Convergefnmf} is not $0$. Then, there exists $\epsilon > 0 $ such that for any $k \in \mathbb{N}$ there are $n_{k} , m_{k} \geq k$ such that
\begin{equation}
\label{Eq:SubseqIntBiggerEpsilon}
\int_{E}|f_{n_{k},m_{k}}-f|d\mu \geq \epsilon. 
\end{equation}
The sequence $f_{k} := f_{n_{k},m_{k}}$ converges $\mu$-almost everywhere to $f$, with $\sup_{k \in \mathbb{N}}|f_{k}|$  integrable. By traditional dominated convergence theorem \citep[VI.9]{doob1953stochastic}, $f_{k} \to f $ in $\mathscr{L}^{1}(E,\mathscr{E},\mu)$. But this contradicts \eqref{Eq:SubseqIntBiggerEpsilon}.  $\blacksquare$

\begin{lemma}
\label{Lemma:IntZmuRiemann}
Let $Z = (Z(x))_{x \in \Rd}$ be a centred second-order stochastic process with locally bounded covariance function. Let $\mu \in \mathscr{M}(\Rd)$. Suppose $Z$ is mean-square continuous outside a $|\mu|$-null set. Let $A \in \BoundedBorel{\Rd}$. For every $n$, let $(I_{j}^{n})_{j \in J_{n}} \subset \Borel{A}$ be a finite partition of $A$ such that $\max_{j \in J_{n}} \diam(I_{j}^{n}) \to 0 $ as $n \to \infty$. Let $x_{j}^{n} \in I_{j}^{n}$ be an arbitrary tag-point for every $(j,n)$. Then the limit in mean-square
\begin{equation}
\label{Eq:IntZmuRiemann}
\int_{A}Z(x)d\mu(x) := \lim_{n \to \infty}\sum_{j \in J_{n}} Z(x_{j}^{n}) \mu( I_{j}^{n} )
\end{equation}
exists and does not depend upon the choice of partitions $(I_{j}^{n})_{j,n}$ or tag-points $(x_{j}^{n})_{j,n}$. In addition, the application $A \mapsto \int_{A}Z(x)d\mu(x)$ is $\sigma$-additive on $\BoundedBorel{\Rd}$.
\end{lemma}

\textbf{Proof of Lemma \ref{Lemma:IntZmuRiemann}:} For a given sequence of Riemann sums, we study the Cauchy gaps for $n,m \in \mathbb{N}$
\small
\begin{equation}
\label{Eq:CauchyGapIntZdmu}
\begin{aligned}
\mathbb{E}&\left( \Big| \sum_{j \in J_{n}}Z(x_{j}^{n}) \mu(I_{j}^{n} ) - \sum_{k \in J_{m}}Z(x_{k}^{m}) \mu(I_{k}^{m} )  \Big|^{2} \right) = \mathbb{E}\left( \Big| \sum_{j \in J_{n}}\sum_{k \in J_{m}}\left( Z(x_{j}^{n}) - Z(x_{k}^{m}) \right) \mu(I_{j}^{n}\cap I_{k}^{m} )  \Big|^{2} \right)   \\
&= \sum_{j,j' \in J_{n}}\sum_{k,k' \in J_{m}} \mathbb{E}\left( (Z(x_{j}^{n}) - Z(x_{k}^{m}))(Z(x_{j'}^{n}) - Z(x_{k'}^{m})  ) \right)\mu(  I_{j}^{n}\cap I_{k}^{m} )\mu(  I_{j'}^{n}\cap I_{k'}^{m} ) \\
&= \int_{A \times A}  \underbrace{\sum_{j,j' \in J_{n}}\sum_{k,k' \in J_{m}}  \left[ C_{Z}( x_{j}^{n} , x_{j'}^{n} ) - C_{Z}( x_{j}^{n} , x_{k'}^{m} ) - C_{Z}( x_{k}^{m} , x_{j'}^{n} ) + C_{Z}( x_{k}^{m} , x_{k'}^{m} ) \right]  \mathbf{1}_{   I_{j}^{n}\cap I_{k}^{m}  \times  I_{j'}^{n}\cap I_{k'}^{m} }}_{f_{n,m} := }  \   d \mu\otimes \mu.
\end{aligned} 
\end{equation}
\normalsize
Let $E \in \Borel{\Rd}$ a set such that $|\mu|(E^{c}) = 0 $ and $Z$ is mean-square continuous over $E$. $C_{Z}$ must thus be continuous over $E\times E$. Therefore, by construction of the Riemann partitions we have $\lim_{n,m \to \infty}f_{n,m}(x,y) = 0$ for every $ (x,y) \in ( E\times E) \cap (A\times A)$. In addition, since $C_{Z}$ is locally bounded and $A$ is bounded, we have
\begin{equation}
\int_{A\times A} \sup_{n,m \in \mathbb{N}} |f_{n,m}| d \left(|\mu|\otimes|\mu|\right)   \leq   4  \sup_{(x,y) \in A \times A} |C_{Z}(x,y)| |\mu|(A)^{2} < \infty.
\end{equation}
By Lemma \ref{Lemma:LDCTdoubleSequence}, \eqref{Eq:CauchyGapIntZdmu} must go to $0$ as $n,m \to \infty$. The sequence of Riemann sums is Cauchy and thus converges to a random variable in $\LtwoOmega$, noted as in \eqref{Eq:IntZmuRiemann}. If we consider another sequences of Riemann sums with partitions $(\tilde{I}_{j}^{n})_{j,n}$ and tag-points $(\tilde{x}_{j}^{n})_{j,n}$, $j \in \tilde{J}_{n}$, then,
\begin{equation}
\label{Eq:DiffGapRiemannIntZmu}
\sum_{j \in J_{n}}Z(x_{j}^{n})\mu(I_{j}^{n}) - \sum_{\tilde{j} \in \tilde{J}_{n}} Z(\tilde{x}_{\tilde{j}}^{n})\mu( \tilde{I}_{\tilde{j}}^{n} ) = \sum_{j \in J_{n}} \sum_{\tilde{j} \in \tilde{J}_{n}} \left( Z(x_{j}^{n}) - Z(\tilde{x}_{\tilde{j}}^{n}) \right) \mu(I_{j}^{n} \cap \tilde{I}_{\tilde{j}}^{n}).
\end{equation}
Applying $\mathbb{E}( \left| \ \cdot \ \right|^{2} )$ to \eqref{Eq:DiffGapRiemannIntZmu}, we can use the same splitting and dominated convergence arguments as in \eqref{Eq:CauchyGapIntZdmu} and conclude that \eqref{Eq:DiffGapRiemannIntZmu} converges to $0$ in $\LtwoOmega$. Thus, the limit does not depend upon the partitions and tag-points.

Finally, using a particular Riemann partition of $A\cup B$ one concludes for $A$ and $B$ bounded
\begin{equation}
\label{Eq:CovIntZmuAxB}
\begin{aligned}
\Cov\left( \int_{A}Z(x)d\mu(x) \ , \ \int_{B} Z(x) d\mu(x) \right) &= \lim_{n \to \infty} \Cov\left( \sum_{j \in J_{n}}Z(x_{j'}^{n})\mu(I_{j'}^{n}\cap A) , \sum_{j' \in J_{n}}Z(x_{j'}^{n})\mu(I_{j'}^{n}\cap B) \right) \\
&=\lim_{n \to \infty} \int_{A\times B}\sum_{j \in J_{n}} \sum_{j' \in J_{n}} C_{Z}( x_{j}^{n} , x_{j'}^{n} ) \mathbf{1}_{I_{j}^{n}\times I_{j'}^{n}} \ d\mu\otimes \mu \\
&= \int_{A\times B} C_{Z} d\mu \otimes \mu,
\end{aligned} 
\end{equation}
which holds from dominated convergence. The application $A \mapsto \int_{A}Z(x)d\mu(x)$ has thus a covariance structure identified with the measure over $\Rd\times\Rd$ given by $D \mapsto \int_{D} C_{Z} d\mu \otimes \mu $, being thus $\sigma$-additive on $\BoundedBorelRd$ (Proposition \ref{Prop:MsigmaAdditive}). $\blacksquare$
 
For the case of an integral over an unbounded set, we define it through growing bounded sets.

\begin{lemma}
\label{Lemma:IntZmuUnbounded}
Let $Z$ and $\mu$ as in Lemma \ref{Lemma:IntZmuRiemann}. Suppose in addition
\begin{equation}
\label{Eq:CzMuOMuIntegrable}
\int_{\Rd\times\Rd} |C_{Z}| \  d|\mu| \otimes |\mu | < \infty.
\end{equation}
Let $A \in \Borel{\Rd}$. Then, for every sequence $(K_{n})_{n \in \mathbb{N}} \subset \BoundedBorel{\Rd}$ with $K_{n} \nearrow \Rd $, the limit in mean-square
\begin{equation}
\label{Eq:DefIntZmuUnbounded}
\int_{A}Z(x)d\mu(x) = \lim_{n \to \infty} \int_{K_{n}\cap A }Z(x)d\mu(x)
\end{equation}
exists and is independent of the growing sequence $(K_{n})_{n}$. In addition, the application $A \mapsto \int_{A}Z(x)d\mu(x)$ is $\sigma$-additive on $\Borel{\Rd}$.
\end{lemma}

\textbf{Proof of Lemma \ref{Lemma:IntZmuUnbounded}:} For a growing sequence $(K_{n})_{n}$, we use the additivity of the integral (Lemma \ref{Lemma:IntZmuRiemann}):
\begin{equation}
\begin{aligned}
\left|\int_{K_{n}}Z(x)d\mu(x) - \int_{K_{m}} Z(x)d\mu(x)\right| &= \left|\int_{ K_{n\vee m} \setminus K_{n \wedge m} }Z(x)d\mu(x)\right|.
\end{aligned}
\end{equation}
Applying $\mathbb{E}( \left| \ \cdot \ \right|^{2} )$, using formula \eqref{Eq:CovIntZmuAxB} we obtain
\begin{equation}
\label{Eq:IntCzmumuDiffKs}
\int_{ K_{n\vee m} \setminus K_{n \wedge m} \times K_{n\vee m} \setminus K_{n \wedge m} } C_{Z} \ d\mu \otimes \mu.
\end{equation}
Since $C_{Z} \in \mathscr{L}^{1}(\RdxRd , |\mu| \otimes |\mu|)$ (condition \eqref{Eq:CzMuOMuIntegrable}), \eqref{Eq:IntCzmumuDiffKs} goes to $0$ as $n,m \to \infty$. Thus, the limit \eqref{Eq:DefIntZmuUnbounded} exists as a limit of a Cauchy sequence in $\LtwoOmega$. Now, if $(\tilde{K}_{n})_{n}$ is another sequence of bounded Borel sets growing to $\Rd$, we have
\begin{equation}
\label{Eq:DiffGapIntZmuUnbounded}
\left| \int_{K_{n}}Z(x)d\mu(x) - \int_{\tilde{K}_{n}}Z(x)d\mu \right| = \left| \int_{K_{n} \bigtriangleup \tilde{K}_{n}} Z(x) d\mu(x)\right|,
\end{equation}
where $K_{n} \bigtriangleup \tilde{K}_{n}$ denotes the symmetric difference between $K_{n}$ and $\tilde{K}_{n}$. The $\mathbb{E}( \left| \ \cdot \ \right|^{2} )$ of \eqref{Eq:DiffGapIntZmuUnbounded} is given by
\begin{equation}
\int_{ K_{n} \bigtriangleup \tilde{K}_{n} \times K_{n} \bigtriangleup \tilde{K}_{n} } C_{Z} \ d\mu \otimes \mu,
\end{equation}
which goes to $0$ as $n \to \infty$ since $C_{Z} \in \mathscr{L}^{1}(\RdxRd , |\mu|\otimes|\mu|)$. Hence, the limit does not depend upon the chosen sequence $(K_{n})_{n}$. Finally, for any $A,B \in \Borel{\Rd}$ we have
\begin{equation}
\begin{aligned}
\Cov\left( \int_{A}Z(x)d\mu(x) \ , \ \int_{B} Z(x) d\mu(x) \right) &= \lim_{n \to \infty} \Cov\left( \int_{A\cap K_{n}} Z(x) d\mu(x) \ , \ \int_{B \cap K_{n}} Z(x)d\mu(x)  \right) \\
&=\lim_{n \to \infty} \int_{A\cap K_{n}\times B\cap K_{n}} C_{Z} \  d\mu\otimes \mu \\
&= \int_{A\times B} C_{Z} d\mu \otimes \mu.
\end{aligned} 
\end{equation}
The application $D \mapsto \int_{D}C_{Z}d\mu$ defines a finite measure and therefore $A \mapsto \int_{A}Z(x) d\mu(x)$ defines a finite m-cov random measure over $\Rd$, being thus $\sigma$-additive over $\Borel{\Rd}$ (Proposition \ref{Prop:MsigmaAdditive}). $\blacksquare$

\textbf{Proof of semi-stochastic Fubini Theorem \ref{Theo:StochasticFubini}:} Condition \eqref{It:Fubini1} implies that the function $\int_{\Rm} \psi( \cdot ,u)d\mu(u) \otimes\int_{\Rm} \psi( \cdot ,v)d\mu(v) $ is in $L^{1}(\RdxRd , |C_{M}|)$\footnote{$\int_{\Rm} |\psi|( \cdot ,u)d|\mu|(u) \otimes\int_{\Rm} |\psi|( \cdot ,v)d|\mu|(v) $ is the Radon-Nikodym derivative of the measure over $\RdxRd$ $$ D \mapsto \int_{D} \int_{\RmxRm } |\psi|(x,u)|\psi|(y,v) d( |\mu| \otimes |\mu| )(u,v)  d|C_{M}|(x,y)$$ with respect to $|C_{M}|$.}. Hence, the integral of  $\int_{\Rm}\psi( \cdot ,u)d\mu(u)$ with respect to $M$ over $\Rd$ (left side of \eqref{Eq:FubiniTheorem}) is a well-defined stochastic integral (condition \eqref{Eq:VarphiIntegrableM}).

For the iterated integral at the right side of \eqref{Eq:FubiniTheorem}, we set
\begin{equation}
Z(u) := \int_{\Rd}\psi(x,u)dM(x).
\end{equation}
The covariance function of $Z$ is
\begin{equation}
\label{Eq:CovCZasIntPsiCM}
C_{Z}(u,v) = \int_{\Rd\times \Rd}\psi(x,u)\psi(y,v) dC_{M}(x,y).
\end{equation}
By condition \eqref{It:Fubini2}, $C_{Z}$ is well-defined and locally bounded, so $Z$ is well-defined as a second-order process. From condition \eqref{It:Fubini1} we conclude $C_{Z} \in \mathscr{L}^{1}(\RmxRm , |\mu|\otimes|\mu|)$ and that $C_{Z}$ is continuous over $E\times E$, being $E \in \Borel{\Rm}$ such that $|\mu|(E^{c}) = 0 $. By Lemmas \ref{Lemma:IntZmuRiemann} and \ref{Lemma:IntZmuUnbounded} the iterated stochastic integral
\begin{equation}
\label{Eq:IntZmuIteratedM}
\int_{\Rm} Z(u) d\mu(u) = \int_{\Rd} \left[ \int_{\Rd} \psi( x,u) dM(x) \right] d\mu(u)
\end{equation}
is well-defined through Riemann sums.

Let us now consider the variance of the difference between the iterated integrals
\small
\begin{equation}
\label{Eq:VarianceDiffFubini}
\begin{aligned}
&\Var\left( \int_{\Rm} \int_{\Rd} \psi( x,u) dM(x) d\mu(u) - \int_{\Rd}\int_{\Rm} \psi( x,u) d\mu(u) dM(x)    \right) =  \Var\left(  \int_{\Rm} \int_{\Rd} \psi( x,u) dM(x) d\mu(u) \right)  \\
&\quad+ \Var\left( \int_{\Rd}  \int_{\Rm} \psi( x,u) d\mu(u) dM(x)  \right) -2\Cov\left(  \int_{\Rm} \int_{\Rd} \psi( x,u) dM(x) d\mu(u)  , \int_{\Rd}\int_{\Rm} \psi( x,u) d\mu(u) dM(x)  \right).
\end{aligned}
\end{equation}
\normalsize
For the variances we have (Eq. \eqref{Eq:CovCZasIntPsiCM} and \eqref{Eq:CovarianceIntegrals})
\begin{equation}
\label{Eq:VarIntMmu}
\Var\left(  \int_{\Rm} \int_{\Rd} \psi( x,u) dM(x) d\mu(u) \right) = \int_{\RmxRm} \int_{\RdxRd} \psi(x,u)\psi(y,v) dC_{M}(x,y) d\mu\otimes \mu (u,v),
\end{equation}
\begin{equation}
\label{Eq:VarIntmuM}
\Var\left( \int_{\Rd}  \int_{\Rm} \psi( x,u) d\mu(u) dM(x) \right) = \int_{\RdxRd} \int_{\Rm} \psi(x,u)d\mu(u)   \int_{\Rm} \psi(y,v) d\mu(v) dC_{M}(x,y).
\end{equation}
For the covariance in \eqref{Eq:VarianceDiffFubini}, we consider a sequence of bounded subsets $(K_{n})_{n}$ growing to $\Rm$ and, for each $n$, a sequence of Riemann partitions of $K_{n}$, say $(I_{j}^{\tilde{n}})_{j \in J_{\tilde{n}},\tilde{n} \in \mathbb{N}}$ and $u_{j}^{\tilde{n}} \in I_{j}^{\tilde{n}}$ tag-points. Then we have
\begin{equation}
\label{Eq:DevelopCovarianceFubiniProof}
\begin{aligned}
\Cov&\left(  \int_{\Rm} \int_{\Rd} \psi( x,u) dM(x) d\mu(u)  , \int_{\Rd}\int_{\Rm} \psi( x,u) d\mu(u) dM(x)  \right) \\
&= \lim_{n \to \infty} \lim_{\tilde{n} \to \infty} \sum_{j \in J_{\tilde{n}}}  \Cov\left( \int_{\Rd}\psi(x,u_{j}^{\tilde{n}}) dM(x) ,   \int_{\Rd}\int_{\Rm} \psi(x,u) d\mu(u) dM(x)    \right) \mu(I_{j}^{\tilde{n}}) \\
&= \lim_{n \to \infty} \lim_{\tilde{n} \to \infty} \sum_{j \in J_{\tilde{n}}} \int_{\RdxRd}  \psi(x,u_{j}^{n}) \int_{\Rm}\psi(x,u) d\mu(u)  dC_{M}(x,y) \mu(I_{j}^{\tilde{n}}) \\
&=\lim_{n \to \infty} \int_{K_{n}} \int_{\RdxRd}  \psi(x,u) \int_{\Rm}\psi(x,u) d\mu(u)  dC_{M}(x,y) d\mu(v) \\
&= \int_{\Rm} \int_{\RdxRd}  \psi(x,u) \int_{\Rm}\psi(x,u) d\mu(u)  dC_{M}(x,y) d\mu(v),
\end{aligned}
\end{equation}
Now, since $C_{M}$ is a measure, from condition \eqref{It:Fubini1} we can use the deterministic Fubini Theorem to argue that all the iterated integrals in \eqref{Eq:VarIntMmu}, \eqref{Eq:VarIntmuM}, \eqref{Eq:DevelopCovarianceFubiniProof} coincide. The variance \eqref{Eq:VarianceDiffFubini} equals thus $0$, showing that both iterated integrals coincide (the random variables involved are all zero-mean). $\blacksquare$

\subsection{Proof of Proposition \ref{Prop:CMprimitiveContinuous}}
\label{App:CMprimitiveContinuousProof}

We consider the following Lemma.

\begin{lemma}
\label{Lemma:PrimitiveMeasureContinuous}
Let $\mu \in \mathscr{M}_{F}(\Rd)$. Then, the function $\vec{x} \mapsto \mu\left( (-\infty , \vec{x}]  \right)$ is continuous outside a Lebesgue measure null set.
\end{lemma}

\textbf{Proof of Lemma \ref{Lemma:PrimitiveMeasureContinuous}:} We first remark that for $d = 1$ this holds immediately since the function $x \mapsto \mu\left( (-\infty , x] \right)$ is càdlàg and therefore it has an at-most countable set of discontinuities (the atoms of $\mu$). For $d > 1$ we proceed as follows. For each $j = 1 , \ldots , d$, consider the positive measure $\mu_{j}$ over $\R$ defined by marginalizing the measure $|\mu|$ over all its components except the $j$-th one:
\begin{equation}
\label{Eq:DefmujProof}
\mu_{j}(A) := |\mu|( \R \times \ldots \times \underbrace{A}_{\hbox{position} \ j} \times \ldots \times \R ), \quad \forall A \in \Borel{\R}.
\end{equation}
Let $F : \Rd \to \R$ be the function $F(\vec{x}) = \mu( (-\infty , \vec{x} ] )$. Then, by additivity
\small
\begin{equation}
\label{Eq:InequalitiesFproof}
\begin{aligned}
| F(\vec{x}+\vec{h}) - F(\vec{x}) | &\leq | F(x_{1} + h_{1} , \ldots , x_{d-1} + h_{d-1} , x_{d} + h_{d}) - F(x_{1} + h_{1} , \ldots , x_{d-1} + h_{d-1} , x_{d} ) |  \\
&\quad+| F(x_{1} + h_{1} , \ldots , x_{d-1} + h_{d-1} , x_{d} ) - F(x_{1} + h_{1} , \ldots , x_{d-1}, x_{d} )  | \\
&\quad+ \ldots \\
&\quad+| F(x_{1} + h_{1} ,  x_{2} , \ldots , x_{d} ) - F(x_{1} , \ldots , x_{d} )  | \\
&= |\mu\big( \ (-\infty , x_{1} + h_{1}] \times \ldots \times(-\infty , x_{d-1} + h_{d-1}] \times  (  x_{d} \wedge ( x_{d} + h_{d} ) ,  x_{d} \vee (x_{d} + h_{d}) ]\ \big)| \\
&\quad + |\mu\big( \ (-\infty , x_{1} + h_{1}]\times\ldots \times  (  x_{d-1}\wedge (x_{d-1} + h_{d-1})  ,    x_{d-1}\vee ( x_{d-1} + h_{d-1}) ] \times (-\infty , x_{d}] \ \big)| \\
&\quad + \ldots \\
&\quad + |\mu\left( \ (  x_{1}\wedge (x_{1} + h_{1})  ,    x_{1}\vee ( x_{1} + h_{1}) ]\times ( -\infty , x_{2} ] \times \ldots \times ( -\infty , x_{d} ] \  \right)  | \\
&\leq \mu_{1}\left( \ \left(x_{1} - |h_{1}| , x_{1} + |h_{1}|\right] \ \right) + \ldots +  \mu_{d}( \ (x_{d} - |h_{d}| , x_{d} + |h_{d}|] \ ).
\end{aligned}
\end{equation} 
\normalsize
Since all the measures $\mu_{j}$ are over $\R$, they all have an at-most countable quantity of atoms. Let $A_{j}$ be the set of atoms of $\mu_{j}$. Note that if $x \notin A_{j}$, then $\mu_{j}( \ (x - |h| , x + |h|] \ ) \to 0 $ as $h \to 0$. Set
\begin{equation}
E := \Big[ \bigcup_{j=1}^{d} \R \times \ldots \times \underbrace{A_{j}}_{\hbox{position } j} \times \ldots \times \R \Big]^{c}.
\end{equation}
Since $E^{c}$ is a finite union of $\ell^{\otimes d}$-null sets, we have $\ell^{\otimes d}(E^{c}) = 0$. From inequality \eqref{Eq:InequalitiesFproof} it follows that $F$ is continuous over $E$. $\blacksquare$

\textbf{Proof of Proposition \ref{Prop:CMprimitiveContinuous}:} From Lemma \ref{Lemma:PrimitiveMeasureContinuous}, there exists $E \in \Borel{\Rd}$ such that $\ell^{\otimes d}(E^{c})=0$ and such that $\vec{u} \mapsto |C_{M}|( \Rd \times (-\infty , \vec{u}]  )$ is continuous over $E$. Using the symmetry of $C_{M}$ and analogue inequalities as in \eqref{Eq:InequalitiesFproof}, we obtain
\small
\begin{equation}
\label{Eq:InequalitiesCmProof}
\begin{aligned}
\left| C_{M}( (-\infty , \vec{u} + \vec{h}_{1} ]\times  (-\infty , \vec{v} + \vec{h}_{2} ] ) - C_{M}\left( (-\infty , \vec{u} ]\times  (-\infty , \vec{v} ] \right) \right| &\leq \left| C_{M}\left( (-\infty , \vec{u}] \ \times \  (-\infty , \vec{v} + \vec{h}_{2} ] \bigtriangleup  (-\infty , \vec{v} ] \right)   \right| \\
&\quad + \left| C_{M}\left( (-\infty , \vec{u} + \vec{h}_{1} ] \bigtriangleup  (-\infty , \vec{u} ]  \ \times \ (-\infty , \vec{v}]  \right)   \right| \\
&\leq |C_{M}|\left( \Rd \ \times \   (-\infty , \vec{u} + \vec{h}_{1} ] \bigtriangleup  (-\infty , \vec{u} ] \right)  \\
&\quad + |C_{M}|\left( \Rd \ \times \   (-\infty , \vec{v} + \vec{h}_{2} ] \bigtriangleup  (-\infty , \vec{v} ] \right), 
\end{aligned}
\end{equation}
\normalsize
By continuity of $\vec{u} \mapsto |C_{M}|( \Rd \times (-\infty , \vec{u}]  )$, if $(\vec{u},\vec{v}) \in E\times E$ then \eqref{Eq:InequalitiesCmProof} goes to $0$ as $(\vec{h}_{1} , \vec{h}_{2} ) \to 0 $. $\blacksquare$

\subsection{Properties of the anti-derivative operator $\mathcal{O}$}
\label{App:PropertiesO}

The derivative property \eqref{Eq:DerOmuMuIntegral} is justified by a convolution argument which we make precise for  $d=1$, the case $d > 1$ is analogous but requiring a more tedious notation. Let $\mu \in \mathscr{M}_{F}(\R)$. We can re-write $\mathcal{O}(\mu)$ as
\begin{equation}
\mathcal{O}(\mu) = \mathbf{1}_{[0, \infty)} \ast \left[ \left( \mathbf{1}_{[0,\infty)} \ast \mu \right) \mathbf{1}_{[0,\infty)} \right] - \mathbf{1}_{(-\infty , 0 )} \ast \left[ \left( \mathbf{1}_{[0,\infty)} \ast \mu \right) \mathbf{1}_{(-\infty , 0 )} \right], 
\end{equation}
where $\ast$ denotes the convolution operation. Using the properties of the convolution with respect to derivatives, that the Dirac measure $\delta$ (at $0$) is its identity element, and that $\frac{d}{dx} \mathbf{1}_{[0,\infty)} = \delta$, one has
\begin{equation}
\begin{aligned}
\frac{d}{dx}\left\{  \frac{d}{dx} \mathcal{O}(\mu) \right\} &= \frac{d}{dx}\left\{   \delta \ast \left[ \left( \mathbf{1}_{[0,\infty)} \ast \mu \right) \mathbf{1}_{[0,\infty)} \right] - (-\delta) \ast \left[ \left( \mathbf{1}_{[0,\infty)} \ast \mu \right) \mathbf{1}_{(-\infty , 0 )} \right] \right\} \\
&=  \frac{d}{dx}\left\{  \left( \mathbf{1}_{[0,\infty)} \ast \mu \right) \mathbf{1}_{[0,\infty)} +  \left( \mathbf{1}_{[0,\infty)} \ast \mu \right) \mathbf{1}_{(-\infty , 0 )}  \right\} \\
&= \frac{d}{dx}\left\{  \mathbf{1}_{[0,\infty)} \ast \mu \right\} \\
&=  \delta \ast \mu = \mu.
\end{aligned}
\end{equation}

For the stochastic case, the derivative relation  \eqref{Eq:OMDerMIntegrals} comes from semi-stochastic Fubini Theorem \ref{Theo:StochasticFubini}: if $\varphi \in \DRd$, we claim that we can exchange the order of integration and conclude

\begin{equation}
\begin{aligned}
\int_{\Rd}\mathcal{O}(M)(\vec{x}) \frac{\partial^{2d} \varphi}{\partial x_{1}^{2} ... \partial x_{d}^{2}}(\vec{x})d\vec{x} &= \int_{\Rd} \int_{0}^{\vec{x}} \int_{\Rd} \mathbf{1}_{(-\infty , \vec{u}]}(s) dM(s) d\vec{u} \frac{\partial^{2d} \varphi}{\partial x_{1}^{2} ... \partial x_{d}^{2}} (\vec{x}) d\vec{x}  \\
&= \int_{\Rd} \int_{\Rd} \int_{0}^{\vec{x}}  \delta_{s}( (-\infty , \vec{u} ] )d\vec{u} \frac{\partial^{2d} \varphi}{\partial x_{1}^{2} ... \partial x_{d}^{2}} (\vec{x}) d\vec{x} dM(s) \\
&= \int_{\Rd} \langle \mathcal{O}( \delta_{s} ) , \frac{\partial^{2d} \varphi}{\partial x_{1}^{2} ... \partial x_{d}^{2}} \rangle dM(s) \\
&= \int_{\Rd}\langle \delta_{s} , \varphi \rangle dM(s) = \int_{\Rd} \varphi(s)dM(s).
\end{aligned} 
\end{equation}
The arguments which verify that the hypotheses in semi-stochastic Fubini Theorem \ref{Theo:StochasticFubini} hold are completely analogous to the case presented in the proof of Theorem \ref{Theo:KLMeasure} (see the proof of Lemma \ref{Lemma:E(M(A)Xj)}). We therefore omit them.

\end{appendix}

\bibliography{mibib}
\bibliographystyle{apacite}

\end{document}